\documentclass[11pt]{amsart}
\input colordvi
\usepackage{amssymb}
\usepackage{amsmath}

\numberwithin{equation}{section}

\newcommand{\pmat}{\begin{pmatrix}}
\newcommand{\epmat}{\end{pmatrix}}

\newcommand{\lk}{\left}

\newcommand{\rk}{\right}
\newcommand{\la}{{\langle}}
\newcommand{\ra}{{\rangle}}

\newcommand{\be} {\begin{enumerate} }
\newcommand{\bcase}{\begin{cases}}
\newcommand{\ecase}{\end{cases}}

\newcommand{\ee} {\end{enumerate} }

\newcommand{\DD}{{\rm I \kern -0.2em D}}
\newcommand{\dd}{{\rm I \kern -0.2em D}}
\newcommand{\NN}{{\rm I \kern -0.2em N}}

\newcommand{\TT}{{\rm I \kern -0.2em T}}

\newcommand{\CV}{{{ \mathcal V }}}

\newcommand{\DEQS}{\begin{eqnarray*}}
\newcommand{\EEQS}{\end{eqnarray*}}
\newcommand{\DEQSZ}{\begin{eqnarray}}
\newcommand{\EEQSZ}{\end{eqnarray}}

\newcommand\del[1]{}
\newcommand\think[1]{}
\newcommand\new[1]{}
\newcommand\zus[1]{}

\newcommand\comc[1]{#1} 
\newcommand\comd[1]{} 
\newcommand\Redd[1]{} 

\def\bdm{\begin{displaymath}}
\def\edm{\end{displaymath}}
\def\bea{\begin{eqnarray}}
\def\eea{\end{eqnarray}}

\newtheorem{theorem}{Theorem}[section]
\newtheorem{lem}[theorem]{Lemma}
\newtheorem{defn}[theorem]{Definition}

\newtheorem{coro}[theorem]{Corollary}

\newtheorem{Rem}{Remark}

\begin{document}

\date{\today}

\title[On a space discretization scheme  for FSHE.]
{On a space discretization scheme  for the Fractional Stochastic
Heat Equations}



\author[Debbi L. ]{Latifa  Debbi}

\address{Department of Mathematics, Chair of Applied Mathematics,
Montanuniversit$\ddot{a}$t, Franz-Josef-Street 18, A-8700 Leoben,
Austria. \& Laboratory of Pure and Applied Mathematics, Faculty of
Sciences, University Ferhat Abbas, El-Maabouda Setif 19000,
Algeria.}

\email{ldebbi@yahoo.fr,  debbi@iecn.u-nancy.fr and
latifa.Debbi@unileoben.ac.at}


\author[Dozzi M.]{Marco Dozzi}

\address{Institut Elie Cartan, Nancy 1 B.P 239, 54506 Vandoeuvre-Les-Nancy
cedex, France.} \email{Dozzi@iecn.u-nancy.fr}

\maketitle \vspace{-1cm} \footnote{This work is supported by The
Austrian Science Foundation grant P26017002 and by the "CNRS Accord
programme" 21378, between Institut Elie Cartan and University Ferhat
Abbas, Setif.}
\begin{abstract}

In this work, we introduce a new discretization to the fractional
Laplacian and use it to elaborate an approximation scheme for
fractional heat equations perturbed by a multiplicative cylindrical
white noise. In particular, we estimate the rate of convergence.



\end{abstract}



\section{Introduction}\label{sec-intro}

In this work, we are interested in the space approximation of the
solutions of  fractional stochastic heat equations. Equations where
leading operators are fractional or more general pseudo-differential
operators are widely used to model complex phenomena. For example,
they are ubiquitous in the study of the quasi geostrophic flow,  the
fast rotating fluids,  the dynamic of the frontogenesis\footnote{The
frontogenesis is the terminology used by atmosphere scientists for
describing the formation in finite time of a discontinuous
temperature front.}  in meteorology,  the diffusions in fractal or
disordered medium, the pollution problems, the mathematical finance
and the transport problems, see e.g. \cite{Bilerandal,
Caffarelli-2009, Caffarelli-Vasseur2010, LeonZas-02,
Pablo-VazquezA-Arxiv, Sug, SugKak, ZasABD-95} and the references
therein. The wellposedness of these equations, in the deterministic
and stochastic cases, has been extensively studied see e.g.
\cite{Bilerandal, BrzezniakDebbi1, Caffarelli-2009,
Caffarelli-Vasseur2010, DebbiDozzi1, Mueller-98, TrumanWu-06}.
Although the numerical approximation of the solution is needed in
applications, the number of  numerical schemes relevant to such
approximations is quite restricted. The main difficulty of the
numerical approximation of fractional equations is related to the
fractional operator. For example, contrarily to second order
differential operators, the fractional operators can not be
discretized by three points.\del{ the fractional operators are
global, hence they can not be discretized by three points.} Using
the classical schemes and as a global operator, all the values on
the grid should be used in every step. Via the integro-differential
representation of the fractional operator, a direct discretization
is based on the discretization of the integrals. This idea has been
used in the numerical study of the deterministic conservation law
driven by fractional power of the Laplacian in \cite{Droniou-09}.
Unfortunately, as it is mentioned in the paper, the convergence of
the scheme elaborated is slow and  leads  to some unreasonable
values. Moreover, in the theoretical study, the author did not give
an explicit form to the discretized operator. The discretization of
the integrals has been already used for the Liouville and Riemann
fractional operators and has yielded  the Gr\"unwald formula, see
for short list \cite{Anh-09, Meerschaert-06, Anh-Tuner-08}.

During the preparation of this work, we found the work of Westphal
\cite{Westphal-73}, on an approach to fractional powers of operators
via fractional differences. The approximation given by Westphal
provides a rigorous mathematical support to a numerical
discretization of fractional operators. Westphal defined the
fractional operator of the infinitesimal generator $A$ of a
semigroup $S(t)$, defined on Banach space X, by:
\begin{equation*}
A^r := s-\lim_{\tau\searrow 0}\tau^{-r}(I-S(\tau))^r, 0<r<1,
\end{equation*}
where
$$(I-S(\tau))^r := \sum_{j=0}^{\infty}
C_j^{-r-1}S(j\tau)=\sum_{j=0}^{\infty}
\frac{\Gamma(j-r)}{\Gamma(j+1)\Gamma(-r)}S(j\tau)$$ and $s-\lim$
means the strong limit, i.e. the limit in the space X of
$\tau^{-r}(I-S(\tau))^r f$, for every $f \in X$. In particular, she
proved:
\begin{equation}\label{Discr-Westphal-Liouville-Fract-Ope}
\lim_{\tau\searrow
0}\tau^{-r}\sum_{j=0}^{[x/\tau]}C_j^{-r-1}f(x-j\tau)= D^rf(x)
\end{equation}
where $ [x/\tau]$ is the integer part of $x/\tau$ and $ D^r $ is the
fractional differential Riemann-Liouville operator defined on $
\mathbb{R}_+$. An intuitive way to discretize the operator $ D^r $
can be obtained by taking $\tau = \frac1n$:
$$ D_n^r f(.) :=
n^{r}\sum_{j=0}^{[nx]}C_j^{-r-1}f(.-\frac{j}{n}).$$

To encounter the difficulties of the direct discretization of the
fractional operator, probabilistic technics have been used. In
particular, in \cite{Woyczynski-Num-approx-fract-05} the authors
used the Monte Carlo method to approximate numerically the solution
of some deterministic fractional partial differential equations,
among them the Burgers equations. In \cite{Schwab-10} the authors
used wavelet techniques to approximate the Kolmogorov  equation
driven by the infinitesimal generator of a Feller process.

Our idea to discretize a fractional operator $ A^\gamma$ is to
discretize first the operator $ A $ then to take the fractional
power of the discrete operator. As far as the authors know, this
idea is new.\del{  The abstract schemes, we know, for parabolic
stochastic equations consider only the global Lipschitz.
Furthermore, for the abstract schemes, operator-space discretization
rather than rate of convergence are not given explicitly.}

In this work, we discretize  the fractional Laplacian, in the way
described above,  and we elaborate a scheme to approximate the
fractional stochastic heat equation. Our aim is also to calculate
explicitly the rate of the convergence and to show its dependence on
the fractional power of the Laplacian. We are also interested in the
critical values of the fractional order which insures the
convergence of the scheme as well.

\vspace{0.5cm} We consider the following fractional stochastic heat
equation:
\begin{equation}\label{Pricipal_Problem_with_x}
\left\{
\begin{array}{rl}
\frac{\partial}{\partial t}u(t,x)&= \frac{\partial^{\alpha}
}{\partial x^{\alpha}}u(t, x)+ g(u(t,x)) \,\frac{\partial^{2}
W}{\partial
t\partial x}(t, x), \,t> 0,x\in (0,1),\\
u(0,x) &= u_{0}(x),\; x\in (0,1),\\
u(t, 0)&=u(t, 1)=0,\; t>0,
\end{array}
\right.
\end{equation}
where $\frac{\partial^{\alpha}}{\partial x^{\alpha}}=
(-\Delta)^{\frac\alpha2}$,  $\alpha>1$ is the fractional power of
the Laplacian. Let us denote by $ A_{\alpha}= A^{\frac{\alpha}{2}}$,
where $A=-\Delta$ is the Laplacian with boundary Dirichlet
conditions, defined on $D(A)=H^{2,2}(0,1)\cap H_0^{1,2}(0,1)$.
$H^{k,p}(0,1)$, for $k \in \mathbb{N}$,  $p \in [1,\infty)$ is the
Sobolev space of order k. The fractional operator
$A_\alpha=A^{\frac{\alpha}{2}}$ is well defined, see e.g. Lemma
2.6.6 in \cite{Pazy} and it is given by the formula (see \cite{Pazy}
pp 72-73):
\begin{equation}\label{Eq-Fract-Pazy}
A_\alpha = (-\Delta)^{\frac\alpha2}: = \frac{\sin \frac{\alpha
\pi}{2} }{\pi}\int_0^\infty z^{\frac{\alpha}{2}-1}A(Iz+A)^{-1}dz.
\end{equation}
The operator $A_\alpha$ is a closed densely defined operator with
domain of definition  given via the complex interpolation of order
$\frac{\alpha}{2}$: $D(A_\alpha)=[H,D(A)]_{\frac{\alpha}{2}}$, see
e.g. \cite{LM-72-i, Triebel-95} and Theorem 4.2 in \cite{Ta}. More
precisely,

$$D(A_\alpha)=D(A^{\alpha/2})=\{v\in L^2(0,1): \sum_{k=1}^\infty
\lambda_k^{\alpha}v_k^2<\infty\},$$ where $v_k=\langle
v,e_k\rangle={\sqrt2}\int_0^1\, v(x)\,\sin{k\pi x}\, dx$ and
$\lambda_k=k^2\pi^2$, $k\in\mathbb{N}$ are  the eigenvalues of the
operator $A$ corresponding to the eigenfunctions:
$\{e_k={\sqrt2}\sin{k\pi \cdot}\}_{ k\in \mathbb{N}}$. The map $ g:
\mathbb{R}\to \mathbb{R}$, is a \comc{bounded } Lipschitz continuous
function on $ \mathbb{R}$. The operator $ g $ is regarded as a
nonlinear operator from $ H= L^{2}(0, 1) $ to $ \mathcal{L}(H) $,
the set of bounded linear operators on $ H $, defined by 
$g(u)(h)=\{(0,1)\ni x \mapsto g(u(x))h(x)\in \mathbb{R}\}$. In other
words, the nonlinear operator $g$ is the Nemytski map associated
with function $g$. For $ v \in H $, $ g(v) $ is given as a
multiplicative operator. From the hypothesis that $ g $ is bounded
we have  $ \Vert g(v)\Vert \leq b_0$, where $b_0= \sup_\mathbb{R}
\vert g(x)\vert$.  $ \{ W(t), t\geq 0\} $ is a cylindrical Wiener
process on the probability space $ (\Omega, \mathcal{F},
\{\mathcal{F}\}_{t\geq 0}, \mathbb{P}) $. The initial condition $
u_{0}$ is a $ L^{2}(0, 1)$-valued $\mathcal{F}_0$-measurable
function. In section \ref{sec-conv-schem}, we will suppose stronger
condition on the diffusion term $g$ and on the initial condition $
u_0$. In particular we will suppose that $ A^\delta g$ is a bounded
Nemytski map for some $ \delta>0$; $ \Vert A^\delta g(v)\Vert \leq
b_\delta$ with $b_\delta= \sup_\mathbb{R} \vert A^\delta g(x)\vert$
and $ u_0$ belongs to a given fractional Sobolev space.

We rewrite the equation (\ref{Pricipal_Problem_with_x}) in the
following form:
\begin{equation}\label{Principal Eq in Heat}
\Big\{
\begin{array}{rl}
du(t)&= -A_{\alpha}u(t)\,dt + g(u(t))\,dW(t), \,t> 0,\\
u(0) &= u_{0}.
\end{array}
\end{equation}

Let us denote by  $ \{ S_{\alpha}(t), t\geq 0 \} $ the semigroup
generated by $ -A_{\alpha} $ and by $H^{\nu, 2}$ the fractional
Sobolev space of order $ \nu$. By a solution of the equation
(\ref{Principal Eq in Heat}), we mean, see e.g. \cite{DaPZa-92}:

\begin{defn} Suppose that $
\alpha>1$. An  $ \mathcal{F}_{t}$-adapted $H^{\theta, 2}$-valued
continuous process $u= \{u(t), t\geq 0\}$ is called a mild solution
of equation (\ref{Principal Eq in Heat}) with initial condition $
u_0 \in H^{\eta, 2}, \eta >0$, iff for some
$p>\frac{2\alpha}{\alpha-1-2\theta}$ where $ \theta <
\min\{\frac{\alpha-1}{2}, \eta\}$
\begin{equation}\label{cond-1}
\mathbb{E}\sup_{t \in [0, T]}|u(t)|_{H^{\theta, 2}}^{\comc{p}} <
\infty, \quad T>0\end{equation}
 and for all $t\geq 0$, a.s. the following identity holds
\begin{equation}\label{Principal Eq in Heat integ form}
u(t)= S_{\alpha}(t)u_{0}+ \int_{0}^{t}S_{\alpha}(t-s)g(u(s))\,dW(s).
\end{equation}
\end{defn}

We introduce the space:

\begin{defn}\label{def-Z_T}
Let $T>0$ and $p\in [1,\infty]$ be fixed and $H$ a Hilbert space. By
$ Z_{T, \theta, p}(H)$ we denote the space of all $H^{\theta,
2}$-valued continuous and $ \mathcal{F}_{t}$-adapted processes $u=
\{u(t), t\in [0, T] \}$ such that
\begin{equation}\label{def-norm-Z_T}
\Vert u\Vert _{T, \theta, p}^{p}:= \mathbb{E}\sup_{t \in [0,
T]}|u(t)|_{H^{\theta, 2}}^{p} = \Vert u\Vert _{L^p(\Omega,
L^\infty(0, T)\times H^{\theta, 2})}< \infty.
\end{equation}
If $\theta =0$, we use shortly the notation $ Z_{T,p}(H)$. \del{When
$X= L^{2}(0, 1)$ we will usually denote the space $ Z_{T, p,
\theta}(L^2(0, 1))$  by $ Z_{T,p, \theta}$ and the norm $\Vert \cdot
\Vert _{T, L^2(0, 1), p}$ by $\Vert \cdot\Vert _{T,p}$.}
\end{defn}

The following result of existence and uniqueness of the solution of
equation (\ref{Principal Eq in Heat}) can be concluded from the
calculus in \cite{BrzezniakDebbi1}:

\begin{theorem}\label{th-main}
Let $\alpha >1$ and let $u_0$ be a $H^{\eta, 2}$-valued
$\mathcal{F}_0$-measurable function such that  $$ \mathbb{E}\vert
u_0\vert^{p}_{H^{\eta, 2}}<\infty $$ for some
$p>\frac{2\alpha}{\alpha-1}$, $ 0\leq \theta <\min
\{\frac{\alpha-1}{2}-\frac\alpha p, \eta\}$ and let $T>0$. Then
there exists a unique mild solution $ u \in Z_{T, \theta, p}$, of
equation (\ref{Principal Eq in Heat}).\del{ Furthermore, $\{u_n(t),
t\geq 0\}$ belongs to the space $C((0, T]; D(A_\alpha^\delta))\cap
C([0, T]; L^2(0, 1))$ provided $ \delta +\frac3{2\alpha}<1$ and $
\delta +p^{-1}<\frac12(1-\frac1\alpha)$.}

\end{theorem}

\del{\begin{theorem}\label{th-main} Let $\alpha >1$ and let $u_0$ be
a $H^{\eta, 2}$-values $\mathcal{F}_0$-measurable function such that
\comc{$$ \mathbb{E}\vert u_0\vert^{p}_{H^{\eta, 2}}<\infty $$} for
some $p>\frac{2\alpha}{\alpha-1-2\theta}$, $ \theta <\min
\{\frac{\alpha-1}{2}, \theta\}$ and let $T>0$. Then there exists a
unique mild solution $ u \in Z_{T, \theta, p}$, of equation
(\ref{Principal Eq in Heat}). Furthermore, \Red{$\{u_n(t), t\geq
0\}$ belongs to the space $C((0, T]; D(A_\alpha^\delta))\cap C([0,
T]; L^2(0, 1))$ provided $ \delta +\frac3{2\alpha}<1$ and $ \delta
+p^{-1}<\frac12(1-\frac1\alpha)$.}
\end{theorem}
} The paper is organized in the following way. In section
\ref{Discretization-Fract-Operator} we describe the discretization
of the fractional operator. In particular, we apply the idea for the
Galerkin approximation and for the finite difference method. In
section \ref{Discretization-Fract-Eq} we elaborate a numerical
scheme to approximate the solution of the fractional stochastic heat
equation (\ref{Principal Eq in Heat}). In section
\ref{sec-prelim-estimations}, we give some preliminary estimations
of the Green functions corresponding to the fractional operator and
to the approximated operator. The section \ref{sec-conv-schem} is
devoted to prove the convergence of the approximated solution to the
solution of the equation (\ref{Principal Eq in Heat}). In the end of
this introduction, let us mention the following references where the
the approximations of certain stochastic partial differential
equations are treated \cite{GongyAlabert-06, Gyongy-Millet-09,
Gyongy-Millet-05, Hausenblas-03, Printems-01}. Let us also mention
that, we take in the whole paper $p, \alpha>1$ and the values of the
constants may change from line to line.


\section{Discretization of the fractional operator}\label{Discretization-Fract-Operator}
Let us first recall notions about the  approximations of the
Laplacian $ -A$, see e.g. \cite{Brenner-Scott}. We consider the
Gelfand triple $(V, H, V')$, where $V\hookrightarrow
H=H'\hookrightarrow V'$ densely, where $V=H_0 ^1(0, 1)$ and $V'$ is
its dual. The operator $A$ defines a coercive bilinear form on $V$
by
$$
a(u,v) := \la A u, v\ra = \la  u', v'\ra, \quad u,v\in V,
$$
where $u', v'$ are the first derivatives of $u$ and $ v$ in the
distribution sense. It is widely accepted that to approximate the
Laplacian $ -A$, it is sufficient to approximate the bilinear form $
a $.

Let $ \CV_n$ be  a finite dimensional subspace of $V$ generated by a
basis $ (f_j)_{j=1}^{n}$. It is easily seen from the formula: $
a(u_n, v_n)= \sum_{i,j=1}^n u_i v_j a(f_i, f_j)$, for $
u_{\CV_n}:=\sum_{j=1}^nu_jf_j , v_{\CV_n}:= \sum_{j=1}^nv_jf_j \in
\CV_n$ that the projection of the bilinear form $ a $ on $\CV_n$,
denoted by $ a_{\CV_n}:= aP_{\CV_n}$ is well defined via the double
series index $ (a(f_i, f_j))_{i,j}$. Using Riesz representation we
can rewrite the  bilinear form $a_{\CV_n}$ as:$$ a_{\CV_n}
(u_{\CV_n}, v_{\CV_n})= \la A_{\CV_n} u_{\CV_n}, v_{\CV_n}\ra=
\sum_{i,j=1}^n u_i v_j \la A_{\CV_n}f_i, f_j\ra,
$$
where $ A_{\CV_n}$ is a positive bounded linear operator on $ \CV_n
$ which is well defined via the stiffness matrix
$\mathbb{A}_{\CV_n}:= (a(f_i, f_j)=\la A_{\CV_n}f_i,
f_j\ra)_{i,j=1}^n$.

Now we define the fractional power of the approximated operator
$A_{\CV_n}$ by the following formula (see \cite{Pazy} pp 72-73):
\DEQSZ\label{form-app} A_{\CV_n}^{\frac\alpha2} =
{\sin\pi{\frac\alpha2}\over \pi} \int_0^{\infty}
z^{{\frac\alpha2}-1}A_{\CV_n}\lk( zI+A_{\CV_n}\rk)^{-1}\, dz \EEQSZ
and the fractional bilinear form
$$ a_{\CV_n, \alpha} (u_{\CV_n}, v_{\CV_n})=
\la A_{\CV_n}^{\frac\alpha2} u_{\CV_n}, v_{\CV_n}\ra.
$$
The fractional stiffness matrix $\mathbb{A}_{\CV_n}^{\frac\alpha2}$
is then given by $ (a_{\CV_n, \alpha}(f_i, f_j)=\la
A_{\CV_n}^{\frac\alpha2}f_i, f_j\ra)_{i,j=1}^n$.

Our idea is to investigate how and on what rate the operator
$A_{\CV_n}^\alpha$ and the bilinear form $ a_{\CV_n, \alpha}$ are
good approximations to the operator $A^\alpha$ respectively the
bilinear form $ a_{\alpha}:= \la A^{\frac\alpha2} u, v\ra$.\del{ if
as much as $ A_{\CV_n}$ is a good approximation of $A$.
\DEQSZ\label{formula-exact} A^\alpha = {\sin\pi\alpha\over \pi}
\int_0^{\infty} z^{-\alpha}A\lk( zI+A\rk)^{-1}\, dz. \EEQSZ
Thus, the approximation $A_h$ of $A$ is given by}

Before going through this calculus, let us apply this method to
calculate the approximation of the fractional operator and of the
stiffness matrix corresponding to the following two methods; the
Galerkin method and the finite difference method.

\subsection{Approximation by Galerkin method.} Let $\CV_n$ be the
subspace generated by the basis $ (e_j)_{j=1}^n$ defined above.
Recall that $e_j(.):= \sqrt 2\sin (j\pi.)$. It is known that the
operator $ A$ and $ A^{\frac\alpha2}$ are diagonizable under the
basis $ (e_j)_{j=1}^\infty$. Consequently the approximating operator
$A_{\CV_n} $ and the stiffness matrix $ \mathbb{A}_{\CV_n}$ are
diagonal with respect to the basis $ (e_j)_{j=1}^n$ and with
corresponding eigenvalues $ (\lambda_i)_{i=1}^n$. Thanks to the
boundness of the approximation operator $A_{\CV_n} $, it is also
easy to define $ A_{\CV_n}$ on the basis as:
$A_{\CV_n}^{\frac\alpha2}e_i= \lambda_i^{\frac\alpha2}e_i,$ see also
Lemma \ref{lem-eigenvalues-fract-and-initial-ops}. \del{$
(\lambda_i^{\frac\alpha2})_{i=1}^n$}

\subsection{Approximation by finite difference method.}
Let $\{ x_i := {i\over n}, i=0,\dots ,n\}$ be the set of grid points
and let $\phi_i, i\in \{1,\dots ,n-1\}$ be a pyramid function, i.e.\
a function, which takes value $1$ at the grid point ${i\over n}$,
vanishes at the other grid points and is linear between the grid
points. The approximating space $\CV_{n-1}$ is then generated by
$\lk\{ \phi_i\rk\} _{i=1} ^{n-1}$. For an implementation reasons, we
will focus on the stiffness matrix. It is well known that the
stiffness matrix $\mathbb{A}_{\CV_{n-1}}$ corresponding to the
finite difference approximation of the operator $ A$ is given by:
$$ a_{ij}= \bcase 2, & i=j,\\
-1, & i=j+1 \mbox{ or } i=j-1, \\ 0, & \mbox{ otherwise.}\ecase
$$
Using the matrix theory it is easy to calculate the fractional power
of the $\mathbb{A}_{\CV_{n-1}}$, denoted by
$\mathbb{A}_{\CV_{n-1}}^{\frac\alpha2}$ via the formula:
\begin{equation}\label{Eq-Fract-Discr-Pazy}
\mathbb{A}_{\CV_{n-1}}^{\frac\alpha2}= \frac{\sin \frac\alpha2\pi
}{\pi}\int_0^\infty
z^{\frac\alpha2-1}\mathbb{A}_{\CV_{n-1}}(Iz+\mathbb{A}_{\CV_{n-1}})^{-1}dz.
\end{equation}
\del{\begin{eqnarray} \la A_{\CV_n}^{\frac\alpha2}e_i, e_j\ra&=&
{\sin\pi{\frac\alpha2}\over \pi} \int_0^{\infty}
z^{{\frac\alpha2}-1}\la A_{\CV_n}\lk( zI+A_{\CV_n}\rk)^{-1}e_i,
e_j\ra \, dz\nonumber \\
&=& {\sin\pi{\frac\alpha2}\over \pi} \lambda_i \int_0^{\infty}
z^{{\frac\alpha2}-1}\lk( z+\lambda_i\rk)^{-1}\la e_i, e_j\ra \,
dz.\nonumber \\
\end{eqnarray}
Hence the fractional stiffness is a diagonal matrix with elements:
$$D_{ii} (\lambda_i)^{\frac\alpha2}{\sin\pi{\frac\alpha2}\over \pi}
\int_0^{+\infty}z^{{\frac\alpha2}-1}\lk( z+\lambda_i\rk)^{-1}dz=
.$$}

\section{Discretization of the fractional stochastic Heat equations}\label{Discretization-Fract-Eq}

Let us first observe that the eigenvalues of
$\mathbb{A}_{\CV_{n-1}}$ are given by $ \lambda_{jn}:=
j^2\pi^2c_{jn}, \; j =1,2,... n-1$ where $ c_{jn}:=
\sin^2(\frac{j\pi}{2n})/(\frac{j\pi}{2n})^2 $,  the corresponding
eigenvectors $ e_j^n = ((e_{jk})_k)$ are given by
\cite{GongyAlabert-06}:
$$ e_{jk}=
\sqrt\frac2n \sin (j\frac{k}{n}\pi).
$$
From Lemma \ref{lem-eigenvalues-fract-and-initial-ops}, it is easy
to see that $  \lambda_{jn}^{\frac\alpha2}\; j =1,2,... n-1$ are the
eigenvalues of $\mathbb{A}_{\CV_{n-1}}^\frac\alpha2$ corresponding
to the eigenvectors $ e_j^n = ((\sqrt\frac2n \sin
(j\frac{k}{n}\pi))_k), \; j =1,2,... n-1$.\del{ By other way,
$\mathbb{A}_{\CV_{n-1}}^{\frac\alpha2}$ can be defined as the
symmetric $n-1-$matrix admits the complete family $ e_j^n \; j
=1,2,... n-1$, as eigenvectors with eigenvalues $
\lambda_{jn}^{\frac\alpha2}\; j =1,2,... n-1$, i.e,
$\mathbb{A}_{\CV_{n-1}}^{\frac\alpha2}$ is equivalent to the
diagonal matrix $\mathbb{A}_{\CV_{n-1}}^{\frac\alpha2} :=
(Dii=\lambda^{\frac\alpha2}_{in}) $, i.e.
$\mathbb{A}_{\CV_{n-1}}^{\frac\alpha2} := P^T D_n^{\frac\alpha2} P$,
where $ P$ is the matrix of eigenvectors.} The semi group
$S_\alpha^n(t)$ generated by $\mathbb{A}_{\CV_{n-1}}^{\frac\alpha2}
$ is given by:
\begin{equation}\label{Eq-Fract-Discr-Sn-alpha}
S_\alpha^n(t)x:= ((\sum_{j=1}^{n-1} G_\alpha^n(t, i, j)x_j)_i),
\end{equation}
where  $G_\alpha^n(t, i, j):=
\sum_{k=1}^{n-1}e^{-t\lambda_{kn}^{\frac\alpha2}}e_{ki}^ne_{kj}^n$
and $ x=(x_j)_{1\leq j\leq n-1}$. Let us define now the operators: $
P_n: L^2(0,1)\rightarrow \mathbb{R}^{n-1} $ and  $ E_n :
\mathbb{R}^{n-1}\rightarrow L^2(0,1)$, called projection
respectively interpolation operators given by the following formula:

For all $ f\in L^2(0,1)$ and for all $ x \in \mathbb{R}^{n-1}$:
\begin{equation}\label{Eq-Def-Pn} P_nf:=
\sum_{k=1}^{n-1}( P_nf)_k e_k^n,
\end{equation}
where
\begin{equation}\label{Eq-Def-Pn}
(P_nf)_k:= \sum_{j=1}^{n-1}
\langle f, e_j \rangle e_j^n(x_k).
\end{equation}
\del{where $ \widetilde{P_n}$ is the projection on $\mathbb{R}^{n-1}$ for
which the components coincide with those given by the projection
$\widehat{P}_n$ on the finite dimensional space in $ L^2$ generated
by $e_1, e_2, ..., e_{n-1}$. i.e. $ \widetilde{P_n}f:= (\langle f,
e_j\rangle)_{j}^T$,  $\widehat{P}_nf:= \sum_{j=1}^{n-1}\langle f,
e_j\rangle e_j$}

\del{\begin{equation}\label{Eq-Def-Pn} (P_nf)_k:= (\sum_{j=1}^{n-1}
\langle f, e_j \rangle e_j^n(k))_k
\end{equation}}

\begin{equation}\label{Eq-Def-En}
E_n x:= \sum_{k=1}^{n-1}\langle x, e^n_k \rangle e_k,
\end{equation}

It is easy to see that the operators $ P_n$ and $ E_n$ satisfy the
properties:
\begin{lem}\label{lem-properties-Pn-En}
\begin{itemize}
\item $ P_n$ and $ E_n$ are bounded  linear operators such that
 $ ||P_n|| \leq 1$ and $ ||E_n|| = 1$.
\item $ P_n E_n =I_n$, where $I_n$ is the identity matrix in
$\mathbb{R}^{n-1}$.
\item $E_nP_n = \widehat{P}_n$, where $\widehat{P}_n$ is the projection
on the finite dimensional space in $ L^2$ generated
by $e_1, e_2, ..., e_{n-1}$. i.e.  $\widehat{P}_nf:= \sum_{j=1}^{n-1}\langle f,
e_j\rangle e_j$.
\item $ P_ne_j= e_j^n \; $, if $j =1,2,... n-1$ and zero if $ j\geq n$.
\item $E_ne_j^n =e_j,  \; j =1,2,... n-1$.
\item $ \{e_j, \; j =1,2,... n-1\} $ are the eigenfunctions
of the operator $ E_n\mathbb{A}_{\CV_{n-1}}^{\frac\alpha2} P_n $
corresponding to the eigenvalues $\lambda_{jn}^{\frac\alpha2}$.
\item The Green function of $ E_n\mathbb{A}_{\CV_{n-1}}^{\frac\alpha2} P_n $, which
is also the kernel of the semigroup: $
E_ne^{-t\mathbb{A}_{\CV_{n-1}}^{\frac\alpha2} }P_n $,  is given by
$$G_\alpha^n(t, x, y):=
\sum_{k=1}^{n-1}e^{-t\lambda_{kn}^{\frac\alpha2}}e_{k}(x)e_{k}(y).$$
\end{itemize}
\end{lem}

\begin{lem}\label{lem-estim-Salpha-discr-sgp}
For $0\leq \delta < \frac14+\frac34\alpha$ and   $ \gamma
>\frac1\alpha-4\frac\delta\alpha$, there exists $ K>0$, such that $ \forall t \in (0, T]$,
\begin{eqnarray}\label{Eq-est-semiGroup-Salpha-Salphan}
|S_\alpha(t)
-E_ne^{-t\mathbb{A}_{\CV_{n-1}}^{\frac\alpha2}}P_n|_{\mathcal{L}(H\rightarrow
D(A^{-\delta}))}&+&\|A^{-\delta}\big(S_\alpha(t)-
E_ne^{-t\mathbb{A}^{\frac\alpha2}_{\mathcal{V}_{n-1}}}P_n\big)\|_{HS}\nonumber
\\ &\leq & K\phi_{\alpha, \delta, \gamma}(t, n),
\end{eqnarray}
where
\begin{equation}\label{Eq-Def-Phi(t, n)}
\phi_{\alpha, \delta, \gamma}(t, n) := \Big\{
\begin{array}{lr}
n^{-\alpha\frac\gamma2-2\delta}t^{-\frac\gamma2}+
n^{-\frac{\alpha}{2}}t^{-\frac{1+\alpha-4\delta}{2\alpha}}, \, \, \, \, \, \,0\leq \delta\leq \frac14\\
n^{-2\delta}+
n^{-\frac{\alpha}{2}}t^{-\frac{1+\alpha-4\delta}{2\alpha}}, \, \, \,
\, \, \, \, \, \, \, \frac14<\delta < \frac14+\frac34\alpha.
\end{array}
\end{equation}

\del{\begin{equation}\label{Eq-Def-Phi(t, n)} \phi_{\alpha, \eta}(t,
n) := \Big\{
\begin{array}{rl}
n^{-\alpha\frac\gamma2}t^{-\frac\gamma2}+
n^{-\frac{\alpha-\alpha\beta}{2}}t^{-\frac{1+\alpha-\beta\alpha}{2\alpha}}, \, \,\eta=0\\
n^{-2\eta}+
n^{-\frac{\alpha-\alpha\beta}{2}}t^{-\frac{1+\alpha-\beta\alpha-4\eta}{2\alpha}},
\, \, \eta\neq 0.
\end{array}
\end{equation}}
\del{\begin{equation}\label{Eq-Def-Phi(t, n)} \phi_{\alpha,
\delta}(t, n) := n^{-2\delta}+
n^{-\frac{\alpha}{2}}t^{-\frac{1+\alpha-4\delta}{2\alpha}}.
\end{equation}}
\end{lem}
Let us make the following convention to write $ \phi_{\alpha,
\delta, \gamma}$ shortly as $ \phi_{\alpha, \delta}$ when $ \gamma $
is not presented, i.e. when $\frac14<\delta <
\frac14+\frac34\alpha$.

\vspace{0.5cm}

Let us now discretize the diffusion term $g$. We denote by $ g_n$
the matrix which is given by the column vectors: $(P_n((g\circ
E_n)e_j))_j, {1\leq j\leq n-1}$, where $ \circ$ is the composition
of the two operators: $ E_n$ and the Nymetsky map $ g $. For $ y\in
\mathbb{R}^{n-1}$, the operator $ (g\circ E_n)y$ acts as the
Nemytski map associated with the function $g(E_ny(.))$, i.e.
$((g\circ E_n)y)h=\{(0,1)\ni x \mapsto g( E_ny(x))h(x)\in
\mathbb{R}\}$. We denote by $ W_n(t)$ the vector $(B_1(t), B_2(t),
..., B_{n-1}(t))$ of independent Brownian motions. We introduce the
following multidimensional stochastic differential equation, where $
u_n(t):= (u_n^k(t))_{1\leq k\leq n-1}$:

\del{. The cylindrical noise $ W(t)= \sum_{j=1}^\infty
B_j(t)\varphi_j(x)$ is discretized by $ W_n(t)= \sum_{j=1}^n
B_j(t)\varphi_j(x)$ $u_n(t):= (u(t, x_k^n))_k=$. We construct the}

\begin{equation}\label{eq-Discret-in-R-n}
\left\{
\begin{array}{rl}
du_n(t)&= -\mathbb{A}_{\CV_{n-1}}^{\frac\alpha2}u_n(t)dt +
g_n(u_n(t))\,dW_n(t), \,t> 0,\\
u_n(0)&:=  P_n(u_{0}).\\
\end{array}
\right.
\end{equation}
For the existence and uniqueness of the solution of the stochastic
differential equation (\ref{eq-Discret-in-R-n}), we refer e.g. to
Theorem 2.3 in \cite{Friedman-75}:
\begin{theorem}\label{Theom.Exist.approx Solu.}
There exists a continuous-$\mathbb{R}^{n-1}$ valued $F_t-$adapted
process $u_{n}= \{u_{n}(t), t\geq 0\}$ solution of the problem
(\ref{eq-Discret-in-R-n}) such that :
\begin{equation}\label{Eq-Solu-integ}
u_n(t)= e^{-t\mathbb{A}_{\CV_{n-1}}^{\frac\alpha2}}u_n(0)+
\int_{0}^{t}e^{-(t-s)\mathbb{A}_{\CV_{n-1}}^{\frac\alpha2}}g_n(u_n(s))dW_n(s),
\;\;\;a.s.
\end{equation}
Furthermore, there exists a constant $ C_{T,n, ||g||}$ such that

\begin{equation}\label{Eq-estimate-solu-discrete} \mathbb{E}
\sup_{t\in [0, T]}|u_{n}(t)|_{\mathbb{R}^{n-1}}^{p} \leq C_{T,n,
||g||}(1+\mathbb{E}|u_0|_{L^{2}}^{p}),  \text{ for all }\;\;\;T>0
\;\;\; \text{and} \;\;\; p\in [1, \infty).
\end{equation}
\end{theorem}

We define the $ L^2-$valued stochastic process $ u^n(t) :=
E_nu_n(t)$. We prove, see Appendix \ref{Appendix-proof-Lemma}, that:

\begin{lem}\label{appendix-E-n-Stochas}
The process $ u^n(t) := E_nu_n(t)$ satisfies the following
stochastic integral equation:
\begin{eqnarray}\label{Eq-Solu-integ-En}
\nonumber u^n(t)&=& E_ne^{-t\mathbb{A}_{\CV_{n-1}}^{\frac\alpha2}
}P_nu_0+
\int_{0}^{t}E_ne^{-(t-s)\mathbb{A}_{\CV_{n-1}}^{\frac\alpha2}}P_ng(u^n(s))dW^n(s), \;\;\;a.s,\\
\end{eqnarray}
where
\begin{eqnarray*}
\int_{0}^{t}E_ne^{-(t-s)\mathbb{A}_{\CV_{n-1}}^{\frac\alpha2}}P_ng(u^n(s))dW^n(s)&:=&
\sum_{j=1}^{n-1}\int_0^t
E_n(e^{-(t-s)\mathbb{A}_{\CV_{n-1}}^{\frac\alpha2}}g_n(u_n(s)))_jdB_j(s)\\
&=&\sum_{j=1}^{n-1}\int_{0}^{t}E_ne^{-(t-s)\mathbb{A}_{\CV_{n-1}}^{\frac\alpha2}}P_ng(u^n(s))e_jdB_j(s)
\end{eqnarray*}
and $(e^{-(t-s)\mathbb{A}_{\CV_{n-1}}^{\frac\alpha2}}g_n(u_n(s)))_j$
is the $j$'s column of the matrix
$e^{-(t-s)\mathbb{A}_{\CV_{n-1}}^{\frac\alpha2}}g_n(u_n(s))$.
Furthermore, for all $T>0 $  and $p\in [1, \infty)$,
$$\mathbb{E} \sup_{t\in [0,
T]}|u^{n}(t)|_{L^2}^{p} \leq C_{T, n,
||g||}(1+\mathbb{E}|u_0|_{L^{2}}^{p}).$$

\end{lem}

In other words, $u^{n}(t)$ satisfies the stochastic partial
differential equation:
\begin{equation*}
\left\{
\begin{array}{rl}
du^{n}(t)&= -E_n\mathbb{A}_{\CV_{n-1}}^{\frac\alpha2}P_nu^{n}(t)\,dt
+
g(u^n(t))dW^n(t), \,t> 0,\\
u^n(0)&=  E_nP_nu_{0}.\\
\end{array}
\right.
\end{equation*}

\del{In the additive case the stochastic term is given by:
\begin{eqnarray*}
\int_{0}^{t}E_ne^{-\mathbb{A}_{\CV_{n-1}}^{\frac\alpha2}(t-s)}dW^n(s):=
\sum_{j=1}^{n-1}\int_0^t E_n
e^{-\mathbb{A}_{\CV_{n-1}}^{\frac\alpha2}(t-s)}P_ne_jdB_j(s),
\end{eqnarray*}}

\section{Preliminary estimates}\label{sec-prelim-estimations}
In this section we give a priori estimations to the Green functions,
$G^n_\alpha$ and $G_\alpha$, corresponding to $
\frac{\partial}{\partial t}- A^{\frac\alpha2}$ respectively $
\frac{\partial}{\partial t}- \mathbb{A}_{\CV_{n-1}}^{\frac\alpha2}$
and to their difference.

\begin{lem}\label{lem-baic-est-G}
For $0\leq \delta < \frac14+\frac34\alpha$ and for all  $ \gamma
>\frac1\alpha-4\frac\delta\alpha$, there exists  $ K>0$, such that $ \forall t \in (0, T]$,

\begin{equation}\label{Eq-est-exp-expn}
\sum_{j=1}^{n-1}\lambda_j^{-2\delta}|e^{-t\lambda_j^{\frac\alpha2}}-
e^{-t\lambda_{jn}^{\frac\alpha2}}|^2 \leq K
n^{-\alpha}t^{-1-\frac1\alpha+4\frac\delta\alpha}
\end{equation}

\begin{equation}\label{Eq-est-exp-infty}
\sum_{j=n}^{\infty}
\lambda_j^{-2\delta}e^{-2t\lambda_j^{\frac\alpha2}} \leq K
\Big(n^{-\alpha\gamma-4\delta}t^{-\gamma
}\mathcal{X}_{[0,\frac14]}(\delta)+
n^{-4\delta}\mathcal{X}_{(\frac14,
\frac14+\frac34\alpha)}(\delta)\Big),
\end{equation}
where $\mathcal{X}_{B} $ is the characteristic function of the set
B: \begin{equation*} \mathcal{X}_{B}(b) := \Big\{
\begin{array}{lr}
1, b\in B\\
0, b\not\in B.\\
\end{array}
\end{equation*}

\end{lem}

\del{\begin{equation}\label{Eq-est-exp-expn}
\sum_{j=1}^{n-1}\lambda_j^{-2\delta}|e^{-t\lambda_j^{\frac\alpha2}}-
e^{-t\lambda_{jn}^{\frac\alpha2}}|^2 \leq K n^{-\alpha
+\alpha\beta}t^{-1-\frac1\alpha+\beta+4\frac\delta\alpha}.
\end{equation}

\begin{equation}\label{Eq-est-exp-infty}
\sum_{j=n}^{\infty}
\lambda_j^{-2\delta}e^{-2t\lambda_j^{\frac\alpha2}} \leq K
\Big(n^{-\alpha\gamma}t^{-\gamma }\delta_0(\delta)+
n^{-4\delta}(1-\delta_0(\delta))\Big),
\end{equation}
where $\delta_0 $ is the Dirac measure in 0.}

\begin{proof}
To get the estimation (\ref{Eq-est-exp-expn}), we use the mean value
theorem. We obtain,

\begin{eqnarray*}
\sum_{j=1}^{n-1}\lambda_j^{-2\delta}|e^{-\lambda_j^\frac\alpha2t}-e^{-\lambda_{jn}^\frac\alpha2t}|^2&=&
\sum_{j=1}^{n-1}\lambda_j^{-2\delta}e^{-2t(\pi j)^\alpha}|1-e^{(1-c_{jn}^\frac\alpha2)(j\pi)^\alpha t}|^2\nonumber\\
&\leq& \sum_{j=1}^{n-1}\lambda_j^{-2\delta}e^{-2t(\pi
j)^\alpha}|1-c_{jn}^\frac\alpha2|^2(j\pi)^{2\alpha} t^2
e^{2(1-c_{jn}^\frac\alpha2)(j\pi)^\alpha t \tau}, \tau \in [0, 1]\nonumber\\
&\leq&
\sum_{j=1}^{n-1}(j\pi)^{-4\delta}|1-c_{jn}^\frac\alpha2|^2(j\pi)^{2\alpha}
t^2
e^{2(-1+(1-c_{jn}^\frac\alpha2))(j\pi)^\alpha t}\nonumber\\
\end{eqnarray*}

We know that:  $ (\frac 2\pi)^\alpha \leq c_{jn}^{\frac\alpha 2}:=
|\sin(\frac{j\pi}{2n})/(\frac{j\pi}{2n})|^\alpha  \leq 1$.
\del{$((-1)^\alpha:=e^{i\alpha\pi})$}Taking $ 1- c_{jn}^{\frac\alpha
2}= O(\frac{j\pi}{2n})^{2\alpha}$, we obtain
\begin{eqnarray*}
\sum_{j=1}^{n-1}\lambda_j^{-2\delta}|e^{-\lambda_j^\frac\alpha2t}-e^{-\lambda_{jn}^\frac\alpha2t}|^2&\leq&
K t^2n^{-4\alpha}\sum_{j=1}^{n-1}j^{6\alpha-4\delta}
e^{-2^{\alpha+1}j^\alpha t}\nonumber\\
&\leq& K t^2n^{-\alpha}\sum_{j=1}^{n-1}j^{3\alpha-4\delta}
e^{-2^{\alpha+1}j^\alpha t}\nonumber\\
&\leq& K t^2n^{-\alpha}\int_{0}^\infty x^{3\alpha-4\delta}
e^{-2^{\alpha+1}x^\alpha t}dx\nonumber\\
&\leq& K
t^{-1-\frac1\alpha+4\frac\delta\alpha}n^{-\alpha}\int_{0}^\infty
y^{3\alpha-4\delta}
e^{-2^{\alpha+1}y^\alpha }dy.\nonumber\\
\end{eqnarray*}
The integral $\int_{0}^\infty y^{3\alpha-4\delta}
e^{-2^{\alpha+1}y^\alpha }dy$  converges provided $ 0\leq \delta <
\frac14+\frac34\alpha$. Hence, there exists $K>0$ such that:

\begin{eqnarray}\label{Eq-Galpha-Galpha-n-partial2}
\sum_{j=1}^{n-1}\lambda_j^{-2\delta}|e^{-\lambda_j^\frac\alpha2t}-e^{-\lambda_{jn}^\frac\alpha2t}|^2\leq K t^{-1-\frac1\alpha+4\frac\delta\alpha}n^{-\alpha}.\nonumber\\
\end{eqnarray}

To get the second estimation (\ref{Eq-est-exp-infty}), let us first
consider the case $ \delta \in[0, \frac14]. $  We use the known
result: for all $ \gamma >0 $, there exists a constant $ k=
k_\gamma$ such that $ e^{-x}\leq K\frac1{x^\gamma}$, we get
\del{(Recall that $ \lambda_j:= (j\pi)^2$)}
\begin{eqnarray}\label{Eq-Galpha-Galpha-n-partial1}
\sum_{j=n}^{\infty}\lambda_j^{-2\delta}e^{-2\lambda_j^{\frac\alpha2}t}
&\leq&
\sum_{j=n}^{\infty}\frac{K}{(2\lambda_j^{\frac\alpha2}t)^\gamma\lambda_j^{2\delta}}
\leq K
n^{-\alpha\gamma-4\delta}t^{-\gamma}\sum_{j=n}^{\infty}(\frac{n}{j})^{\alpha\gamma+4\delta}.
\end{eqnarray}
But  $ \sum_{j=n}^{\infty}(\frac{n}{j})^{\alpha\gamma+4\delta} \leq
K\int_1^{+\infty}x^{-\alpha\gamma-4\delta}dx <\infty,$ provided
$\frac1\alpha-4\frac\delta\alpha<\gamma$. Hence

\begin{eqnarray}\label{Eq-Galpha-Galpha-n-partial10}
\sum_{j=n}^{\infty}\lambda_j^{-2\delta}e^{-2\lambda_j^{\frac\alpha2}t}
&\leq& K n^{-\alpha\gamma-4\delta}t^{-\gamma}.
\end{eqnarray}
For $ \delta \in (\frac14, \frac14+\frac34\alpha)$, we use the
inequality:
$\sum_{j=n}^{\infty}\lambda_j^{-2\delta}e^{-2\lambda_j^{\frac\alpha2}t}
\leq \sum_{j=n}^{\infty}\lambda_j^{-2\delta} $, than we arguing as
above and using the condition $ \delta >\frac14$, we obtain
\begin{eqnarray}\label{Eq-Galpha-Galpha-n-partial1F}
\sum_{j=n}^{\infty}\lambda_j^{-2\delta}e^{-2\lambda_j^{\frac\alpha2}t}
&\leq& \sum_{j=n}^{\infty}\lambda_j^{-2\delta} \leq K
n^{-4\delta}\sum_{j=n}^{\infty}(\frac{n}{j})^{4\delta}\leq K
n^{-4\delta}.
\end{eqnarray}

\end{proof}

\begin{lem}\label{lem-approx-Galpha-Galpha-n}
For $0\leq \delta < \frac14+\frac34\alpha$  and  $ \gamma
>\frac1\alpha-4\frac\delta\alpha$, there exists  $ K>0$, such that $ \forall t \in (0, T]$, we have
\begin{equation}
\sup_{x\in[0,1]}|A_x^{-\delta}\big(G_\alpha (t, x, .)- G_\alpha^n
(t, x, .)\big)|_{H} \leq K\phi_{\alpha, \delta, \gamma }(t, n),
\end{equation}
where $ \phi_{\alpha, \delta, \gamma}(t, n) $ is given by
(\ref{Eq-Def-Phi(t, n)}).
\end{lem}

\begin{proof}
Using the definitions of the functions $G_\alpha (t, x, y)$ and $G_\alpha^n (t, x, y)$ and the fact that the orthonormal basis $ (e_k)_{k\geq1} \subset L^\infty(0,1)$, we get
\begin{eqnarray}\label{Eq-G-Gnstep1}
|A^{-\delta}_x(G_\alpha (t, x, .)&-&G_\alpha^n (t, x,
.))|_{H}\nonumber \\
&=&
|\sum_{j=1}^{n-1}(e^{-\lambda_j^{\frac\alpha2}t}-e^{-\lambda_{jn}^{\frac\alpha2}t})
A^{-\delta}_xe_j(x)e_j+
\sum_{j=n}^{\infty}e^{-\lambda_j^{\frac\alpha2}t}
A^{-\delta}_xe_j(x)e_j|_{H}\nonumber \\
&\leq &
|\sum_{j=1}^{n-1}(e^{-\lambda_j^{\frac\alpha2}t}-e^{-\lambda_{jn}^{\frac\alpha2}t})
A^{-\delta}_xe_j(x)e_j|_{H}+
|\sum_{j=n}^{\infty}e^{-\lambda_j^{\frac\alpha2}t}
A^{-\delta}_xe_j(x)e_j|_{H}\nonumber \\
&\leq &
|\sum_{j=1}^{n-1}(e^{-\lambda_j^{\frac\alpha2}t}-e^{-\lambda_{jn}^{\frac\alpha2}t})
\lambda_j^{-\delta}e_j(x)e_j|_{H}+
|\sum_{j=n}^{\infty}e^{-\lambda_j^{\frac\alpha2}t}
\lambda_j^{-\delta}e_j(x)e_j|_{H}\nonumber \\
& \leq &
\big(\sum_{j=1}^{n-1}(e^{-\lambda_j^{\frac\alpha2}t}-e^{-\lambda_{jn}^{\frac\alpha2}t})^2
\lambda_j^{-2\delta}(e_j(x))^2 \big)^\frac12+
\big(\sum_{j=n}^{\infty}e^{-2\lambda_j^{\frac\alpha2}t}
\lambda_j^{-2\delta}(e_j(x))^2\big)^\frac12 \nonumber \\
& \leq &
\big(\sum_{j=1}^{n-1}\lambda_j^{-2\delta}(e^{-\lambda_j^{\frac\alpha2}t}-e^{-\lambda_{jn}^{\frac\alpha2}t})^2
\big)^\frac12+
\big(\sum_{j=n}^{\infty}\lambda_j^{-2\delta}e^{-2\lambda_j^{\frac\alpha2}t}\big)^\frac12. \nonumber \\
\end{eqnarray}

Replacing (\ref{Eq-est-exp-infty}) and (\ref{Eq-est-exp-expn}) in
(\ref{Eq-G-Gnstep1}), we get the result.

\end{proof}
As a consequence of Lemma \ref{lem-approx-Galpha-Galpha-n}, we obtain:
\begin{coro}\label{Coro-approx-Galpha-Galpha-n}
Under the same conditions in  Lemma
\ref{lem-approx-Galpha-Galpha-n},
\begin{equation}\label{G-Gn-HtimesH}
|A^{-\delta}(G_\alpha (t, ., .)-G_\alpha^n (t, ., .))|_{H\times H}
\leq K \phi_{\alpha, \delta, \gamma}(t, n).
\end{equation}
\end{coro}

\subsection*{Proof of Lemma \ref{lem-estim-Salpha-discr-sgp}} Let $
f\in H$. The semigroups $ S_\alpha $ and
$E_ne^{-t\mathbb{A}_{\CV_{n-1}}^{\frac\alpha2}}P_n $ are acting on
an element $f $ via their Green functions in the following:
\begin{equation*}
(S_\alpha (t)f)(x)=  \int_0^1G_\alpha(t,x ,y)f(y)dy
\end{equation*}
and
\begin{equation*}
(E_ne^{-t\mathbb{A}_{\CV_{n-1}}^{\frac\alpha2}}P_nf)(x)=
\int_0^1G^n_\alpha(t,x ,y)f(y)dy.
\end{equation*}
Applying the H\"older inequality, we get:

\begin{eqnarray*}
|A^{-\delta}(S_\alpha
(t)-E_ne^{-t\mathbb{A}_{\CV_{n-1}}^{\frac\alpha2}}P_n)
f|_{H}&=&\big[\int_0^1|A^{-\delta}((S_\alpha
(t)-E_ne^{-t\mathbb{A}_{\CV_{n-1}}^{\frac\alpha2}}P_n)
f)(x)|^2dx\big]^{\frac 12}
\nonumber \\
&=&\big[\int_0^1|\int_0^1A_x^{-\delta}(G_\alpha(t, x, y)-G_\alpha^n
(t, x, y)) f(y)dy|^2dx\big]^{\frac12}
\nonumber \\
&\leq& \big[\int_0^1dx(\int_0^1|A_x^{-\delta}(G_\alpha (t, x,
y)-G_\alpha^n (t, x, y))|^2dy)\big]^{\frac 12} |f|_H
\nonumber \\
&\leq& \big[\int_0^1|A_x^{-\delta}(G_\alpha (t, x, .)-G_\alpha^n (t,
x, .))|_H^2dx\big]^{\frac 12} |f|_H
\nonumber \\
&\leq& |A_x^{-\delta}(G_\alpha (t, ., .)-G_\alpha^n (t, .,
.)|_{H\times H} |f|_H.
\nonumber \\
\end{eqnarray*}
Using (\ref{G-Gn-HtimesH}), we get the result. For the second
estimation, we have by a direct application of the definitions of
Hilbert-Schmidt norm and the properties of the semigroups $
S_\alpha(t)$ and $E_ne^{-t\mathbb{A}_{\CV_{n-1}}^{\frac\alpha2}
}P_n$,

\begin{eqnarray}
\|A^{-\delta}\big( S_\alpha(t)-
E_ne^{-t\mathbb{A}^{\frac\alpha2}_{\mathcal{V}_{n-1}}}P_n\big)\|_{HS}^2&:=&
\sum_{j=1}^\infty |A^{-\delta}\big(S_\alpha(t)-
E_ne^{-t\mathbb{A}^{\frac\alpha2}_{\mathcal{V}_{n-1}}}P_n\big)e_j|_{H}^2\nonumber\\
&\leq& \sum_{j=1}^{n-1}
\lambda_j^{-2\delta}|e^{-t\lambda_j^\frac\alpha2}-
e^{-t\lambda_{jn}^\frac\alpha2}|^2|e_j|_{H}^2+ \sum_{j=n}^{\infty}
\lambda_j^{-2\delta}e^{-2t\lambda_j^\frac\alpha2}|e_j|_{H}^2\nonumber\\
&\leq& \sum_{j=1}^{n-1}
\lambda_j^{-2\delta}|e^{-t\lambda_j^\frac\alpha2}-
e^{-t\lambda_{jn}^\frac\alpha2}|^2+ \sum_{j=n}^{\infty}
\lambda_j^{-2\delta}e^{-2t\lambda_j^\frac\alpha2}.\nonumber\\
\end{eqnarray}
Thanks to the basic inequality $ (a+b)^\frac12\leq
a^\frac12+b^\frac12$ for $a, b\geq 0$ and by the formula
\eqref{Eq-est-exp-expn} and \eqref{Eq-est-exp-infty}, we get the
result.

\del{\begin{eqnarray}\label{Eq-est-Salpha-Salphan-priori}
\|S_\alpha(t)-
E_ne^{-t\mathbb{A}^{\frac\alpha2}_{\mathcal{V}_{n-1}}}P_n\|_{HS}&\leq&
K\big(n^{-\alpha+\alpha\beta}t^{-1-\frac1\alpha+\beta}+n^{-\alpha\gamma}t^{-\gamma}
\big)^\frac12\nonumber\\
&\leq& K\phi_{\alpha, 0}(t, n)\nonumber\\
\end{eqnarray}}

\section{Convergence of the scheme}\label{sec-conv-schem}
Now we are ready to give the main result of this work.
\begin{theorem}\label{Theom.Approx}
For  $ \alpha>1$,
\begin{equation}\label{eq-cond-eta}
 \frac14+\frac\alpha4< \eta<
\frac14+3\frac\alpha4,
\end{equation}
\begin{equation}\label{eq-cond-delta}
\max\{\frac12, \frac14+\frac\alpha8\}<\delta<  \frac14+\frac34\alpha
\end{equation}
and
 \begin{equation}\label{eq-cond-p}
 p > \max\{\frac{2 \alpha }{\alpha-2}, \frac{ \alpha
}{2\delta-1}, \frac{2 \alpha }{8\delta-\alpha-2}\} \end{equation}
assume that:
\begin{itemize}
\item $(H_g)$ the diffusion term is the Nemytski operator defined by the Lipschitz function  $ g \in D(A^{\delta}) $ and such that $ b_\delta:=\sup_{x\in
\mathbb{R}}\vert A^\delta g(x)\vert<\infty$.
\item $(H_{u_0})$ the initial condition $ u_0$ is an
$D(A^{\eta})-$valued $ L^p$ random variable i.e. $ u_0 \in
L^p(\Omega, D(A^{\eta}))$.
\end{itemize}
Then $u^{n}= \{u^{n}(t), t\in [0,T] \}$ converges to $u:= \{u(t),
t\in [0,T]\}$ in the space $ Z_{T,p}(L^2(0, 1))$. Furthermore, there
exists a constant $ K>0$, such that, \del{\begin{eqnarray}
\|u(t)-u^n(t)\|_{Z_{T,p}(L^2(0, 1)}:= \Big(\mathbb{E}\big[\sup_{[0,
T]}|u(t)-u^n(t)|_H^p\big]\Big)^\frac1p &\leq&
K_{T,|u_0|_{D(A^{\eta})}, b_\delta}n^{-\frac\alpha2},\nonumber \\
\end{eqnarray}}
\begin{eqnarray}\label{Eq-mean-rate-estimation}
\|u(t)-u^n(t)\|_{Z_{T,p}(L^2(0, 1))}&\leq&
K_{T,|u_0|_{D(A^{\eta})}, b_\delta}n^{-\xi},\nonumber \\
\end{eqnarray}
where $ \xi $ is given by
\begin{eqnarray}\label{Eq-xi-1} \xi= \min\{
\frac\alpha2, \; 2\delta\}.
\end{eqnarray}
In particular, for $ 1<\alpha\leq 2$, the rate of convergence $
\xi:= \frac\alpha2.$
\end{theorem}

\begin{Rem}
Let us remark that the rate of convergence is independent of the
regularity of the diffusion term when the dissipation order $
\alpha$ is less than the Laplacian dissipation. In this case, it is
enough to take $ g \in H^1.$
\end{Rem}
\begin{theorem}\label{Theom.Approx-2}
Assume that $ \alpha > 2$, $ p > \frac{2 \alpha }{\alpha-2} $, $g$
and $ u_0 $ satisfy respectively $(H_g)$ with $ \delta =0$ and
$(H_{u_0})$, with $\eta $ satisfying \eqref{eq-cond-eta}. Then
$u^{n}= \{u^{n}(t), t\in [0,T] \}$ converges to $u:= \{u(t), t\in
[0,T]\}$ in the space $ Z_{T,p}(L^2(0, 1))$ and
\begin{eqnarray}\label{Eq-mean-rate-estimation-second}
\|u(t)-u^n(t)\|_{Z_{T,p}(L^2(0, 1))}&\leq&
K_{T,|u_0|_{D(A^{\eta})}, b_0}n^{-(\frac{\alpha}4-\frac12-\frac\alpha{2p})}.\nonumber \\
\end{eqnarray}
\del{ where
\begin{eqnarray}\label{Eq-xi-2} \xi:= \min\{
\frac\alpha2-1-\frac\alpha p, \frac\alpha4-\frac\alpha{2p}, 2\eta\}.
\end{eqnarray}}

\end{theorem}
First let us introduce some lemmata which we will use in the proof
of the convergence.

\begin{lem}\label{lem-Cumm} The operator $ A^{-\delta}$ commutes with $ S_\alpha(t)$ and with
$E_ne^{-\mathbb{A}_{\CV_{n-1}}^{\frac\alpha2}t}P_n $, for all $t
\geq 0$.
\end{lem}

\begin{proof}
For the proof see the Appendix B. \del{\appendix{Proof of
\ref{lem-Cumm}}.}
\end{proof}

\begin{lem}\label{lem-HS-AN}
Suppose that $z\in L^q(0,1)$ for  $q \in [2,\infty]$. Let $Z$ denote
the multiplication operator by $z$. Then, for  $\beta < 1-
\frac{1}{\alpha}$, there exists a constant $ K>0$, such that
\begin{eqnarray}
\int_{0}^\infty s^{-\beta}\Vert
E_ne^{-s\mathbb{A}^\frac\alpha2_{\mathcal{V}_{n-1}}}P_nZ\Vert
_{HS}^{2}\,ds &\leq& K \big(\sum_{k=1}^{+\infty} k^{-\alpha
(1-\beta)}\big)\vert z \vert_{L^q}^2
<\infty.\nonumber \\
\end{eqnarray}
\end{lem}

\begin{proof}
Let us first estimate the term $ \Vert E_ne^{-s\mathbb{A}^\frac\alpha2_{\mathcal{V}_{n-1}}}P_nZ\Vert _{HS}^{2}$. Using Lemma \ref{lem-properties-Pn-En}, we have

\begin{eqnarray*}
\Vert E_ne^{-s\mathbb{A}^\frac\alpha2_{\mathcal{V}_{n-1}}}P_n Z\Vert
_{HS}^{2}&=&\sum_{k=1}^\infty\vert
E_ne^{-s\mathbb{A}^\frac\alpha2_{\mathcal{V}_{n-1}}}P_nZe_k\vert^2_{L^2}=
\sum_{k=1}^\infty\vert Z
E_ne^{-s\mathbb{A}^\frac\alpha2_{\mathcal{V}_{n-1}}}
P_ne_k|^2_{L^2}\\&=&\sum_{k=1}^{n-1}\vert Z
E_ne^{-s\mathbb{A}^\frac\alpha2_{\mathcal{V}_{n-1}}} e^n_k|^2_{L^2}
=\sum_{k}^{n-1}\vert z e^{-\lambda^\frac\alpha2_{kn} s}
E_ne^n_k|^2_{L^2}\\
&=& \sum_{k}^{n-1}e^{-\lambda^\frac\alpha2_{kn} s}\vert z e_k|^2_{L^2}.\\
\end{eqnarray*}
Let us observe that  by the H\"older inequality, $\vert z
e_k\vert_{L^2}\leq \vert z \vert_{L^q} \vert e_k\vert_{L^{r}}$,
where $\frac1{r}+\frac1{q}=\frac12$. Moreover, since $\vert
e_k\vert_{L^2}=1$ and $\vert e_k\vert_{L^\infty}=2^{1/2}$ it follows
by applying the H\"older inequality  that $\vert
e_k\vert_{L^{r}}\leq 2^{1/q}$. Let us recall that $
\lambda^\frac\alpha2_{kn}:= (c_{kn}(k\pi)^2)^{\frac\alpha2}$, where
$ c_{kn}:= \frac{\sin^2(\frac{k\pi}{2n})}{(\frac{k\pi}{2n})^2}$, we
get,
\begin{eqnarray*}
\Vert E_ne^{-s\mathbb{A}^\frac\alpha2_{\mathcal{V}_{n-1}}}P_n Z\Vert _{HS}^{2}&\leq& 2^{\frac{2}{q}}\vert z
\vert_{L^q}^2 \sum_{k=1}^{n-1}e^{-\lambda^\frac\alpha2_{kn} s}\\
&\leq& 2^{\frac{2}{q}}\vert z
\vert_{L^q}^2 \sum_{k=1}^{n-1}e^{-(|c_{kn}|k\pi)^\alpha s}.
\end{eqnarray*}
Therefore and thanks to the fact that: $(\frac 2\pi)^\alpha \leq
|c_{kn}|^{\frac\alpha 2}=
|\sin(\frac{k\pi}{2n})/(\frac{k\pi}{2n})|^\alpha  \leq1 $,
\begin{eqnarray*} \int_{0}^{\infty}s^{-\beta}\Vert
E_ne^{s\mathbb{A}^\frac\alpha2_{\mathcal{V}_{n-1}}}P_nZ\Vert_{HS}^{2}\,ds & \leq  & 2^{\frac2q}\vert z
\vert_{L^q}^2 \int_{0}^{\infty}s^{-\beta} \sum_{k=1}^{n-1}e^{-(c_{kn})^\frac\alpha2(k\pi)^\alpha s}\,ds \nonumber\\
&=&  2^{\comc{\frac2q}}\vert z \vert_{L^q}^2\sum_{k=1}^{n-1}
(c_{kn})^{\frac12(\alpha \beta-\alpha)}(k\pi)^{\alpha
\beta-\alpha}\int_{0}^{\infty} \tau^{-\beta}
e^{-\tau}d\tau \nonumber\\
&\leq& K\vert z \vert_{L^q}^2\sum_{k=1}^{+\infty} (k\pi)^{\alpha
\beta-\alpha}\int_{0}^{\infty} \tau^{-\beta}
e^{-\tau}d\tau \nonumber\\
&\leq& K\vert z \vert_{L^q}^2\sum_{k=1}^{+\infty}
k^{-\alpha (1-\beta)}.\nonumber\\
\end{eqnarray*}
Since by our assumptions $ \alpha (1-\beta)> 1$ the series on the
RHS above is convergent and the result follows.
\end{proof}

\del{\begin{lem}
We first consider the case when $g$ is bounded, i.e. $g$ satisfies
the condition (\ref{eqn-lin_growth}) with $b_1=0$. Let $u \in
Z_{T,p}$. By Lemma \ref{lem-HS} we have

\begin{eqnarray}
\int_{0}^{t}\Vert S_{\alpha}(t-s)g(u(s))\Vert _{HS}^{2}\,ds &\leq&
\int_{0}^{t}\Vert S_{\alpha}(t-s)\Vert _{HS}^{2}\Vert g(u(s))\Vert_{L^\infty}^{2}\,ds \nonumber\\
&\leq& b_0^{2}\int_{0}^{t}\Vert S_{\alpha}(s)\Vert _{HS}^{2}\,ds.
\label{ineq:HS}
\end{eqnarray}
Since $ \alpha > 1 $ and so $ 0 < 1- \frac{1}{\alpha} $ by Corollary
\ref{coro-estim-1/2g} we infer that for    all  $t>0$, the integral
on the RHS of the inequality (\ref{ineq:HS}) is convergent. Hence
the It\^o stochastic integral $\int_0^t S_{\alpha}(t-s)g(u(s))\,
dW(s)$ exists and thus the process $G(u)$ is well defined. We shall
prove that in fact  $ Gu\in Z_{T,p} $.
\end{lem}}

The following Lemma is a special case of Lemma 2.7 from
\cite{BrzezniakDebbi1}:
\begin{lem}\label{lem R gamma}
Provided that $ \nu > p^{-1}$ the operator $ R_{\nu}: L^{p}(0, T;
L^{2}(0, 1)) \rightarrow C([0, T]; L^{2}(0, 1)))$ given by
$$
R_{\nu}h(t)= \int_{0}^{t}(t-s)^{\nu-1}S_{\alpha}(t-s)h(s)\,ds, \; h
\in L^{p}(0, T; L^{2}(0, 1))
$$
is well defined, linear and bounded. Moreover, there exists a
constant $C_{p,\nu}>0$ such that for all $h \in L^{p}(0, T; L^{2}(0,
1))$
\begin{eqnarray}\label{ineq:Rg}
|R_{\nu}h|_{C([0, T]; L^2(0, 1))}&\leq&
K_{p,\nu}T^{\nu-\frac{1}{p}}|h|_{L^{p}(0, T; L^{2}(0, 1))}.
\end{eqnarray}
\end{lem}

\del{Using a similar proof of Lemma 2.7 from \cite{BrzezniakDebbi1}
and the inequalities (\ref{Eq-est-semiGroup-Salpha-Salphan})and
(\ref{Eq-Def-Phi(t, n)}), we prove:}

\begin{lem}\label{lem G gamma} Let $ \frac14<\delta<\frac14+\frac34\alpha$ and $p>\max\{1, \frac{2\alpha}{\alpha-1+4\delta}\}$,  then there
exists $\nu$, satisfying  $ \max\{0, \frac12+
\frac1{2\alpha}-\frac{2\delta}{\alpha}+\frac1p\} <\nu<1$, such that
the operator $ G_{n, \nu}: L^{p}(0, T; D(A^\delta)) \rightarrow
C([0, T]; L^{2}))$ given by
$$
G_{n, \nu}h(t)= \int_{0}^{t}(t-s)^{\nu-1}\big[S_{\alpha}(t-s)-
E_ne^{-(t-s)\mathbb{A}^{\frac\alpha2}_{\mathcal{V}_{n-1}}}P_n\big]h(s)\,ds,
\; h \in L^{p}(0, T; D(A^{\delta}))
$$
is well defined, linear and bounded. Moreover, there exists a
constant $K_{p, T}>0$ such that for all $h \in L^{p}(0, T;
D(A^{\delta}))$
\begin{eqnarray}\label{ineq:Rg}
|G_{n, \nu}h|_{C([0, T]; L^2)}&\leq& K_{T, p}(n^{-2\delta}+
n^{-\frac{\alpha}{2}})|h|_{L^{p}(0, T; D(A^\delta))}.
\end{eqnarray}
\end{lem}

\begin{proof} Let us fix $h \in L^{p}(0, T; D(A^\delta))$. Then for
$t \in (0,T)$ and thanks to Lemma \ref{lem-Cumm}, we have
\begin{eqnarray*}
|G_{n, \nu}h(t)|_{L^2}&=&
|\int_{0}^{t}(t-s)^{\nu-1}\big[S_{\alpha}(t-s)-
E_ne^{-(t-s)\mathbb{A}^{\frac\alpha2}_{\mathcal{V}_{n-1}}}P_n\big]A^{-\delta}A^{\delta}h(s)\,ds|_{L^2}\nonumber\\
&\leq& \int_{0}^{t}(t-s)^{\nu-1}|A^{-\delta}\big[S_{\alpha}(t-s)-
E_ne^{-(t-s)\mathbb{A}^{\frac\alpha2}_{\mathcal{V}_{n-1}}}P_n\big]A^{\delta}h(s)|_{L^2}\,ds\nonumber\\
&\leq& \int_{0}^{t}(t-s)^{\nu-1}\Vert S_{\alpha}(t-s)-
E_ne^{-(t-s)\mathbb{A}^{\frac\alpha2}_{\mathcal{V}_{n-1}}}P_n\Vert_{\mathcal{L}(H\rightarrow D(A^{-\delta}))} |A^{\delta}h(s)|_{L^2}\,ds.\nonumber\\
\end{eqnarray*}
From Lemma \ref{lem-estim-Salpha-discr-sgp},
((\ref{Eq-est-semiGroup-Salpha-Salphan}) and (\ref{Eq-Def-Phi(t,
n)})) and by applying  the H\"older inequality, we get

\begin{eqnarray*} |G_{n, \nu}h(t)|_{L^2}&\leq&
K\int_{0}^{t}(t-s)^{\nu-1}\phi_{\alpha, \delta}(t-s, n)|A^{\delta}h(s)|_{L^2}\,ds\nonumber\\
&\leq& K\big(\int_{0}^{T}(s^{\nu-1}\phi_{\alpha, \delta}(s,
n))^{\frac{p}{p-1}}\,ds\big)^{\frac{p-1}{p}}
\big(\int_{0}^{T}|A^{\delta}h(s)|_{L^2}^{p}\,ds\big)^{\frac{1}{p}}.\nonumber\\
\end{eqnarray*}
Thanks to the basic inequality: $ (x+y)^\theta\leq c_\theta
(x^\theta+ y^\theta)$, for $ \theta >1$ and $ x, y\geq 0,$ we have

\begin{eqnarray}\label{Conv-Int-sPower}
\int_{0}^{T}(s^{\nu-1}\phi_{\alpha, \delta}(s,
n))^{\frac{p}{p-1}}\,ds&=& \int_{0}^{T}(s^{\nu-1}(n^{-2\delta}+
n^{-\frac{\alpha}{2}}s^{-\frac{1+\alpha-4\delta}{2\alpha}}))^{\frac{p}{p-1}}\,ds
\nonumber\\
&\leq &c_p(n^{-2\delta\frac{p}{p-1}}T^\nu
+n^{-\frac{\alpha}{2}\frac{p}{p-1}}\int_0^T
s^{(\nu-1-\frac{1+\alpha-4\delta}{2\alpha})\frac{p}{p-1}}\,ds).
\nonumber\\
\end{eqnarray}
The last integral in the RHS of \eqref{Conv-Int-sPower} converges
provided  $\nu>
\frac12+\frac1{2\alpha}-\frac{2\delta}{\alpha}+\frac1p$. Hence
\del{When $ p>\frac{2\alpha}{2\alpha-1}$, $ \frac1\alpha<\gamma <
2\frac{p-1}{p}$ there exists $\nu> \frac\gamma2+\frac1p$, such that
the first integral (\ref{Conv-Int-sPower}) converges. Furthermore,
if  $ \beta > \frac2p+\frac1\alpha-1$ then we can chose $\nu$
satisfying also the inequality  $\nu>
\frac{1+\alpha-\alpha\beta}{2\alpha}+\frac1p$ and hence the second
integral in (\ref{Conv-Int-sPower}) converges also. Hence for $\nu>
\max\{\frac{1+\alpha-\alpha\beta}{2\alpha}+\frac1p,
\frac\gamma2+\frac1p\}$, we have:}
\begin{eqnarray*}
\big(\int_{0}^{T}(s^{\nu-1}\phi_{\alpha, \delta}(s,
n))^{\frac{p}{p-1}}\,ds\big)^{\frac{p-1}{p}}\leq c_{p,
T}(n^{-2\delta}+ n^{-\frac{\alpha}{2}}).
\nonumber\\
\end{eqnarray*}
The choice of $\nu$ such that $ \frac12+
\frac1{2\alpha}-\frac{2\delta}{\alpha}+\frac1p<\nu<1$ is possible
thanks to the condition $p>\frac{2\alpha}{\alpha-1+4\delta}$.
Finally, we have for all $ t \geq 0 $
\begin{eqnarray*}
|G_{\nu}h(t)|_{L^2}&\leq& K_{T, p}(n^{-2\delta}+
n^{-\frac{\alpha}{2}})|h|_{L^{p}(0, T; D(A^\delta))}. \nonumber
\end{eqnarray*}

\subsection*{Proof of Theorem \ref{Theom.Approx}}

Let
\begin{equation}\label{eq-Mn-2}
M_n:= ||u - u^n||^p_{T, p} = \mathbb{E}\sup_{t \in [0,
T]}|u(t)-u^n(t)|_{H}^{p}.
\end{equation}
We have from equations (\ref{Principal Eq in Heat integ form}) and (\ref{Eq-Solu-integ-En}),
\begin{equation*}
|u(t)-u^n(t)|_{H}^{p} \leq c_p(A(t)+ B(t))
\end{equation*}
where
\begin{eqnarray}\label{eq-A(t)-1}
A(t)&:= &|S_\alpha (t)u_0-E_ne^{-t\mathbb{A}^{\frac\alpha2}_{\mathcal{V}_{n-1}}}P_nu_0|_{H}^{p} \nonumber \\
\end{eqnarray}
\begin{eqnarray}\label{eq-B(t)-1}
B(t)&:=& |\int_0^tS_\alpha(t-s)g(u(s))dW(s)-\int_{0}^{t}
E_ne^{-(t-s)\mathbb{A}^{\frac\alpha2}_{\mathcal{V}_{n-1}}}P_ng(u^n(s))dW^n(s)|_{H}^{p} \nonumber\\
\end{eqnarray}

\subsubsection*{Estimation of $ A(t)$:} Using Lemmata  \ref{lem-Cumm} and \ref{lem-estim-Salpha-discr-sgp}, ((\ref{Eq-est-semiGroup-Salpha-Salphan}) and
(\ref{Eq-Def-Phi(t, n)})), We obtain,
\begin{eqnarray*}
A(t)&:= &|(S_\alpha (t)-E_ne^{-t\mathbb{A}^{\frac\alpha2}_{\mathcal{V}_{n-1}}}P_n)u_0|_{H}^p\nonumber \\
&= &|A^{-\eta}(S_\alpha (t)-E_ne^{-t\mathbb{A}^{\frac\alpha2}_{\mathcal{V}_{n-1}}}P_n)A^{\eta}u_0|_{H}^p\nonumber \\
&\leq&|(S_\alpha
(t) -E_ne^{-t\mathbb{A}^{\frac\alpha2}_{\mathcal{V}_{n-1}}}P_n|^p_{\mathcal{L}(H\rightarrow D(A^{-\eta}))}|A^{\eta}u_0|_{H}^{p}\nonumber \\
&\leq& |u_0|_{D(A^{\eta})}^{p}\phi_{\alpha, \eta}^p(t, n) \nonumber
\\ &\leq& |u_0|_{D(A^{\eta})}^{p}\big(n^{-2\eta}+
n^{-\frac{\alpha}{2}}t^{-\frac{1+\alpha-4\eta}{2\alpha}}\big)^p. \nonumber \\
\end{eqnarray*}
Thanks to the condition: $\eta>\frac14+\frac\alpha4$, the power
$-\frac{1+\alpha-4\eta}{2\alpha}$ is positive, hence
\begin{eqnarray*}
A(t) &\leq& K_T|u_0|_{D(A^{\eta})}^{p}\big(n^{-2\eta}+
n^{-\frac{\alpha}{2}}\big)^p. \nonumber \\
\end{eqnarray*}
Consequently:
\begin{eqnarray}\label{Eq-EsupA}
\mathbb{E}\big[\sup_{[0, T]}A(t)\big] &\leq&
K_T\mathbb{E}|u_0|_{D(A^{\eta})}^{p}\big(n^{-2\eta}+
n^{-\frac{\alpha}{2}}\big)^p \nonumber \\
&\leq& K_T\mathbb{E}|u_0|_{D(A^{\eta})}^{p}\big(n^{-2p\eta}+
n^{-\frac{p\alpha}{2}}\big).\nonumber \\
\end{eqnarray}

\subsubsection*{Estimation of $ B(t)$:}
\del{\begin{eqnarray}\label{Eq-first-Est-B(t)} B(t)&:=&
|\int_0^tS_\alpha(t-s)g(u(s))dW(s)-\int_{0}^{t}
E_ne^{-(t-s)\mathbb{A}^{\frac\alpha2}_{\mathcal{V}_{n-1}}}P_ng(u^n(s))dW^n(s)|_{H}^{p}\nonumber \\
&\leq & c_p\big\{|\sum_{j=1}^{n-1}\int_0^t\big[S_\alpha(t-s)g(u(s))-
E_ne^{-\mathbb{A}_{\CV_{n-1}}^{\frac\alpha2}(t-s)}P_ng(u^n(s))\big]e_jdB_j(s)|_{H}^{p}\nonumber \\
&+ &
|\sum_{j=n}^{\infty}\int_0^tS_\alpha(t-s)g(u(s))e_jdB_j(s)|_{H}^{p}\big\}\nonumber \\
&\leq & c_p\big\{|\sum_{j=1}^{n-1}\int_0^t\big(S_\alpha(t-s)g(u(s))-
E_ne^{-\mathbb{A}_{\CV_{n-1}}^{\frac\alpha2}(t-s)}P_n\big)g(u^n(s))e_jdB_j(s)|_{H}^{p}\nonumber \\
&+ &|\sum_{j=1}^{n-1}\int_0^tS_\alpha(t-s)\big(g(u(s))-
g(u^n(s))\big)e_jdB_j(s)|_{H}^{p} \nonumber \\ &+&
|\sum_{j=n}^{\infty}\int_0^tS_\alpha(t-s)g(u(s))e_jdB_j(s)|_{H}^{p}\big\}\nonumber \\
\end{eqnarray}}

\del{\begin{eqnarray}\label{Eq-first-Est-C(B)-0} B(t)&:=&
|\int_0^tS_\alpha(t-s)g(u(s))dW(s)-\int_{0}^{t}
E_ne^{-(t-s)\mathbb{A}^{\frac\alpha2}_{\mathcal{V}_{n-1}}}P_ng(u^n(s))dW^n(s)|_{H}^{p}.\nonumber \\
\end{eqnarray}}
Let us first introduce the transformations $ \overbrace{.}^{n}$
defined on the  set of Nemytsky maps N, such that for $ h\in N$:
\begin{equation}
\overbrace{h(u)}^ne_j:= \Big\{
\begin{array}{rl}
h(u)e_j\,\, j<n,  \\
0,\,\, j\geq n
\end{array}
\end{equation}
Then we write the second stochastic integral in RHS of
(\ref{eq-B(t)-1}), as
\begin{eqnarray}\label{Eq-beatiful-writting}
\int_{0}^{t}E_ne^{-(t-s)\mathbb{A}^{\frac\alpha2}_{\mathcal{V}_{n-1}}}P_ng(u^n(s))dW^n(s)&:=&\sum_{j=1}^{n-1}
\int_0^tE_ne^{-\mathbb{A}_{\CV_{n-1}}^{\frac\alpha2}(t-s)}P_ng(u^n(s))e_jdB_j(s)\nonumber \\
&=& \int_{0}^{t}E_ne^{-(t-s)\mathbb{A}^{\frac\alpha2}_{\mathcal{V}_{n-1}}}P_n\overbrace{g(u^n(s))}^ndW(s).\nonumber \\
\end{eqnarray}
Using the factorization method see e.g. \cite{BrzezniakDebbi1} and
\cite{DaPZa-92}, we can again rewrite the integrals in
(\ref{eq-B(t)-1}) for $ 0<\nu<1$ as:
\begin{eqnarray*}
\int_0^tS_\alpha(t-s)g(u(s))dW(s)&=&
\int_0^t(t-s)^{\nu-1}S_\alpha(t-s)Y(s)ds
\end{eqnarray*}
and
\begin{eqnarray*}
\int_{0}^{t}E_ne^{-(t-s)\mathbb{A}^{\frac\alpha2}_{\mathcal{V}_{n-1}}}P_n\overbrace{g(u^n(s))}^ndW(s)&=&
\int_0^t(t-s)^{\nu-1}E_ne^{-(t-s)\mathbb{A}^{\frac\alpha2}_{\mathcal{V}_{n-1}}}P_nY_n(s)ds,
\end{eqnarray*}
where
\begin{eqnarray}\label{Eq-Def-Y}
Y(s):=\int_0^s(s-r)^{-\nu}S_\alpha(s-r)g(u(r))dW(r)
\end{eqnarray}
and
\begin{eqnarray}\label{Eq-Def-Yn}
Y_n(s):=\int_0^s(s-r)^{-\nu}E_ne^{-A_n^\alpha(s-r)}P_n\overbrace{g(u^n(r))}^ndW(r)
\end{eqnarray}
Consequently,
\begin{eqnarray}\label{Eq-ReRew-Est-C(B)}
B(t)&=&|\int_0^t(t-s)^{\nu-1}S_\alpha(t-s)Y(s)ds-
\int_0^t(t-s)^{\nu-1}E_ne^{-(t-s)\mathbb{A}^{\frac\alpha2}_{\mathcal{V}_{n-1}}}P_nY_n(s)ds|_{H}^p\nonumber \\
&\leq & c_p(|\int_0^t(t-s)^{\nu-1}S_\alpha(t-s)\big(Y(s)-Y_n(s)\big)ds|_{H}^p\nonumber \\
&+&|\int_0^t(t-s)^{\nu-1}\Big[S_\alpha(t-s)-
E_ne^{-(t-s)\mathbb{A}^{\frac\alpha2}_{\mathcal{V}_{n-1}}}P_n\Big]
Y_n(s)ds|_{H}^p)\nonumber \\
\end{eqnarray}
Using Lemmata \ref{lem R gamma} and \ref{lem G gamma} and taking $$
 \max\{p^{-1}, \frac12+
\frac1{2\alpha}-\frac{2\delta}{\alpha}+\frac1p\}<\nu<1$$ we deduce
that
\begin{eqnarray}\label{Eq-ReRew-Est-C(B)-2}
\mathbb{E}\Big[\sup_{[0,T]}B(t)\Big]&\leq & c_p\mathbb{E}\sup_{[0,T]}\Big(|\int_0^t(t-s)^{\nu-1}S_\alpha(t-s)\big(Y(s)-Y_n(s)\big)ds|_{H}^p\nonumber \\
&+&|\int_0^t(t-s)^{\nu-1}\Big[S_\alpha(t-s)-
E_ne^{-(t-s)\mathbb{A}^{\frac\alpha2}_{\mathcal{V}_{n-1}}}P_n\Big]
Y_n(s)ds|_{H}^p\Big)\nonumber \\
&\leq&    K_{T, p}\Big[\mathbb{E}|Y-Y_n|^p_{L^{p}(0, T; L^{2}(0,
1))}+ (n^{-2\delta}+
n^{-\frac{\alpha}{2}})^p\mathbb{E}|Y_n|^p_{L^{p}(0, T; D(A^{\delta}))}\Big]\nonumber \\
\end{eqnarray}

\subsubsection*{\bf Calculation of $\mathbb{E}|Y_n|^p_{L^{p}(0, T; D(A^\delta))}$}
By the Burkholder's inequality and  Lemma \ref{lem-HS-AN}, there
exists a constant $C_p>0$ such that

\begin{eqnarray}\label{Inq. Y(s) in Lp(0,T; L2) Ex. case}
\int_{0}^{T}\!\!\!\!\!\!&{}&\!\!\!\!\mathbb{E}\big|Y_n(s)
\big|_{D(A^{\delta})}^{p}\,ds \leq
C_p\int_{0}^{T}\!\!\!\mathbb{E}\Big(\int_{0}^{s}(s-r)^{-2\nu}\Vert
A^{\delta}E_ne^{-(s-r)\mathbb{A}^{\frac\alpha2}_{\mathcal{V}_{n-1}}}P_ng(u(r))\Vert
_{HS}^{2}dr
\Big)^{\frac{p}{2}}\,ds \nonumber\\
&\leq &\!\!\!
K_p\int_{0}^{T}\!\!\!\mathbb{E}\Big(\int_{0}^{s}(s-r)^{-2\nu}\Vert
A^{\delta}E_ne^{-(s-r)\mathbb{A}^{\frac\alpha2}_{\mathcal{V}_{n-1}}}P_nA^{-\delta}A^{\delta}g(u(r))\Vert
_{HS}^{2}dr
\Big)^{\frac{p}{2}}\,ds \nonumber\\
\!\!\!&\leq &\!\!\!
K_p\int_{0}^{T}\!\!\!\mathbb{E}\Big(\int_{0}^{s}(s-r)^{-2\nu}\Vert
E_ne^{-(s-r)\mathbb{A}^{\frac\alpha2}_{\mathcal{V}_{n-1}}}P_n\Vert_{HS}^{2}\sup_{r\in[0,
T]}\vert A^{\delta}g(u(r))\vert_{L^\infty(H)}^2dr
\Big)^{\frac{p}{2}}\,ds \nonumber\\
\!\!\! &\leq &\!\!\!
K_pb_\delta^{p}\int_{0}^{T}\!\!\!\Big(\int_{0}^{s}r^{-2\nu}\Vert
E_ne^{-r\mathbb{A}^{\frac\alpha2}_{\mathcal{V}_{n-1}}}P_n\Vert
_{HS}^{2}dr
\Big)^{\frac{p}{2}}\,ds \nonumber\\
\!\!\!&\leq&\!\!\!K_pb_\delta^{p}T\!\Big(\sum_{k=1}^{+\infty}
k^{-\alpha (1-2\nu)}
\Big)^{\frac{p}{2}}.\nonumber\\
\end{eqnarray}
Since  $ p>\max\{\frac\alpha{2\delta-1}, \frac{2\alpha}{\alpha-1}\}$
and $\delta>\frac12$, we can chose $ \nu$, such that
$$\max\{p^{-1}, \frac1p+\frac{1}{2}+ \frac{1}{2 \alpha}- \frac{2\delta}{\alpha}\} <
\nu < \frac{1}{2}- \frac{1}{2 \alpha} $$ what implies that
$\alpha(1-2\nu)>1$ we infer that the last term in \eqref{Inq. Y(s)
in Lp(0,T; L2) Ex. case} is finite.

\subsubsection*{\bf Calculation of $\mathbb{E}|Y-Y_n|^p_{L^{p}(0, T; L^{2}(0, 1))}$}
Using the formula \eqref{Eq-Def-Y} and \eqref{Eq-Def-Yn}, we have:
\begin{eqnarray*}
|Y(s)-Y_n(s)|^p_{H}&\leq& c_p\Big(
\Big|\int_{0}^{s}(s-r)^{-\nu}\big( S_\alpha(s-r)- E_ne^{-(s-r)\mathbb{A}^{\frac\alpha2}_{\mathcal{V}_{n-1}}}P_n\big)g(u(r))dW(r)\Big|_{H}^p \nonumber \\
&+&
\Big|\int_{0}^{s}(s-r)^{-\nu} E_ne^{-(s-r)\mathbb{A}^{\frac\alpha2}_{\mathcal{V}_{n-1}}}P_n\big(g(u(r))-\overbrace{g(u^n(r))}^n\big)dW(r)
\Big|_{H}^p\Big).\nonumber \\
\end{eqnarray*}
By Burkholder's inequality, we have
\begin{eqnarray}\label{est-first-Y-Yn}
\mathbb{E}|Y(s)-Y_n(s)|^p_{H}&\leq& c_p\Big(
\mathbb{E}\Big(\int_{0}^{s}(s-r)^{-2\nu}\|\big( S_\alpha(s-r)- E_ne^{-(s-r)\mathbb{A}^{\frac\alpha2}_{\mathcal{V}_{n-1}}}P_n\big)g(u(r))\|_{HS}^2dr\Big)^\frac p2 \nonumber \\
&+& \mathbb{E}\Big(\int_{0}^{s}(s-r)^{-2\nu}
\|E_ne^{-(s-r)\mathbb{A}^{\frac\alpha2}_{\mathcal{V}_{n-1}}}P_n\big(g(u(r))-\overbrace{g(u^n(r))}^n\big)\|^2_{HS}dr
\Big)^\frac p2 \Big).\nonumber \\
\end{eqnarray}Using the well known functional inequality $ \|AB\|_{HS}\leq
\|A\|_{HS}|B|_{\mathcal{L}(H)}$, we estimate the first term in the
RHS of \eqref{est-first-Y-Yn} as follow:

\begin{eqnarray}\label{est-first-term-in-Y-Yn}
\int_{0}^{s}(s&-&r)^{-2\nu}\|\big( S_\alpha(s-r)- E_ne^{-(s-r)\mathbb{A}^{\frac\alpha2}_{\mathcal{V}_{n-1}}}P_n\big)g(u(r))\|_{HS}^2dr \nonumber \\
&\leq& \int_{0}^{s}(s-r)^{-2\nu}\|A^{-\delta}\big( S_\alpha(s-r)- E_ne^{-(s-r)\mathbb{A}^{\frac\alpha2}_{\mathcal{V}_{n-1}}}P_n\big)A^\delta g(u(r))\|_{HS}^2dr \nonumber \\
&\leq& \int_{0}^{s}(s-r)^{-2\nu}\|A^{-\delta}\big( S_\alpha(s-r)-
E_ne^{-(s-r)\mathbb{A}^{\frac\alpha2}_{\mathcal{V}_{n-1}}}P_n\big)\|_{HS}^2
|A^{\delta}g(u(r))|_{\mathcal{L}(H)}^2dr\nonumber \\
&\leq& b_\delta^2\int_{0}^{s}(s-r)^{-2\nu}\|
A^{-\delta}\big(S_\alpha(s-r)-
E_ne^{-(s-r)\mathbb{A}^{\frac\alpha2}_{\mathcal{V}_{n-1}}}P_n\big)\|_{HS}^2
dr\nonumber \\
\end{eqnarray}

Thanks to \eqref{Eq-est-semiGroup-Salpha-Salphan}
\eqref{Eq-Def-Phi(t, n)}, we have

\begin{eqnarray*}
\int_{0}^{s}(s&-&r)^{-2\nu}\| A^{-\delta}\big(S_\alpha(s-r)-
E_ne^{-(s-r)\mathbb{A}^{\frac\alpha2}_{\mathcal{V}_{n-1}}}P_n\big)\|_{HS}^2
dr\nonumber \\
&\leq&K \int_{0}^{s}(s-r)^{-2\nu}\phi^2_{\alpha, \delta}(s-r, n)dr\nonumber\\
&\leq&K\Big(n^{-4\delta}\int_{0}^{s}r^{-2\nu}dr+ n^{-\alpha}\int_{0}^{s}r^{-2\nu-1-\frac1\alpha+\frac{4\delta}{\alpha}}dr\Big).\nonumber \\
\end{eqnarray*}
The integral $\int_{0}^{s}r^{-2\nu}dr$ is finite thanks to the
condition $ \nu< \frac12-\frac1{2\alpha}<\frac12$ and the integral
$\int_{0}^{s}r^{-2\nu-1-\frac1\alpha+\frac{4\delta}{\alpha}}dr$
converges thanks to the conditions $
\nu<\min\{\frac{2\delta}{\alpha}-\frac1{2\alpha}, \frac{1}{2}-
\frac{1}{2 \alpha}\}$. Hence, we take the parameter $\nu$ which
satisfies the following inequalities:
\begin{eqnarray}\label{Eq-nu-final-cond}
\max\{p^{-1}, \frac1p+\frac{1}{2}+ \frac{1}{2 \alpha}-
\frac{2\delta}{\alpha}\} < \nu <
\min\{\frac{2\delta}{\alpha}-\frac1{2\alpha}, \frac{1}{2}-
\frac{1}{2 \alpha}\}.
\end{eqnarray}
The parameter $ \nu$ exists thanks to the conditions: $\delta
>\frac14+\frac\alpha8$ and $ p>\max\{\frac{2\alpha}{8\delta-2-\alpha}, \frac\alpha{2\delta-1}\}$. Hence
\begin{eqnarray}\label{est-int-term-HS}
\int_{0}^{s}(s&-&r)^{-2\nu}\| S_\alpha(s-r)-
E_ne^{-(s-r)\mathbb{A}^{\frac\alpha2}_{\mathcal{V}_{n-1}}}P_n\|_{HS}^2
dr\nonumber \\
&\leq&K_T\Big(n^{-4\delta}+n^{-\alpha}\Big)\nonumber \\
\end{eqnarray}
By replacing \eqref{est-int-term-HS}in
\eqref{est-first-term-in-Y-Yn}, we get
\begin{eqnarray*}
\int_{0}^{s}(s&-&r)^{-2\nu}\|\big( S_\alpha(s-r)- E_ne^{-(s-r)\mathbb{A}^{\frac\alpha2}_{\mathcal{V}_{n-1}}}P_n\big)g(u(r))\|_{HS}^2dr \nonumber \\
&\leq&
b_\delta^2K_T\Big(n^{-4\delta}+n^{-\alpha}\Big).\nonumber\\
\end{eqnarray*}
Hence,
\begin{eqnarray}\label{est-first-term-in-Y-Yn-final-corr}
\mathbb{E}\Big(\int_{0}^{s}(s&-&r)^{-2\nu}\|\big( S_\alpha(s-r)- E_ne^{-(s-r)\mathbb{A}^{\frac\alpha2}_{\mathcal{V}_{n-1}}}
P_n\big)g(u(r))\|_{HS}^2dr \Big)^\frac p2\nonumber \\
&\leq&
b_\delta^pK_T\Big(n^{-4\delta}+n^{-\alpha}\Big)^\frac p2\nonumber\\
&\leq&
b_\delta^pK_T\Big(n^{-2p\delta}+n^{-p\frac\alpha2}\Big)\nonumber\\
\end{eqnarray}

Now we estimate the second term in \eqref{est-first-Y-Yn}. Arguing
as in the proof of Lemma \ref{lem-estim-Salpha-discr-sgp}, we have
\begin{eqnarray*}
\nonumber \|E_ne^{-(s-r)\mathbb{A}^{\frac\alpha2}_{\mathcal{V}_{n-1}}}&P_n&\big(g(u(r))-\overbrace{g(u^n(r))}^n\big)\|^2_{HS}\\
&=&\|\big(g(u(r))-\overbrace{g(u^n(r))}^n\big)E_ne^{-(s-r)\mathbb{A}^{\frac\alpha2}_{\mathcal{V}_{n-1}}}P_n\|^2_{HS} \nonumber \\
&=&\sum_{j=1}^{\infty}
|\big(g(u(r))-\overbrace{g(u^n(r))}^n\big)E_ne^{-(s-r)\mathbb{A}^{\frac\alpha2}_{\mathcal{V}_{n-1}}}P_ne_j|^2_{H}. \nonumber \\
\end{eqnarray*}
From Lemma \ref{lem-properties-Pn-En}, we have
\begin{eqnarray*}
E_ne^{-(s-r)\mathbb{A}^{\frac\alpha2}_{\mathcal{V}_{n-1}}}P_ne_j=\Big\{
\begin{array}{lr}
e^{-(s-r)\lambda_{jn}^\frac\alpha2}e_j, \,\,\, j\in \{1,
\cdot\cdot\cdot,
n-1\}\\
0, \,\,\,\,\,\,\,\,\,\,\,\,\,\,\,\,\,\,\,\,\,\,\,\,\,\,\,\,\,\,j\geq
n,
\end{array}
\end{eqnarray*}
using the definition of $\overbrace{g(u^n(r))}^n$, the Lipschitz
property of $g$ and Lemma \ref{lem-properties-Pn-En}, we get
\begin{eqnarray*}
\nonumber \|E_ne^{-(s-r)\mathbb{A}^{\frac\alpha2}_{\mathcal{V}_{n-1}}}&P_n&\big(g(u(r))-\overbrace{g(u^n(r))}^n\big)\|^2_{HS}\\
&=&\sum_{j=1}^{n-1}e^{-2(s-r)\lambda_{jn}^\frac\alpha2}
|\big(g(u(r))-g(u^n(r))\big)e_j|^2_{H} \nonumber \\
&\leq&\sum_{j=1}^{n-1}e^{-2(s-r)\lambda_{jn}^\frac\alpha2}
|g(u(r))-g(u^n(r))|^2_{\mathcal{L}(H)} \nonumber \\
&\leq&\sum_{j=1}^{n-1}e^{-2(s-r)\lambda_{jn}^\frac\alpha2}
|u(r)-u^n(r)|^2_{H}. \nonumber \\
&\leq&
\sup_{r\in[0, s]}|u(r)-u^n(r)|^2_{H}\sum_{j=1}^{n-1}e^{-2(s-r)\lambda_{jn}^\frac\alpha2}. \nonumber \\
\end{eqnarray*}
Hence,
\begin{eqnarray}\label{est-second-term-Y-Yn}
\mathbb{E}\Big(\int_{0}^{s}(s&-&r)^{-2\nu}
\|E_ne^{-(s-r)\mathbb{A}^{\frac\alpha2}_{\mathcal{V}_{n-1}}}P_n\big(g(u(r))-\overbrace{g(u^n(r))}^n\big)\|^2_{HS}dr
\Big)^\frac p2 \Big) \nonumber \\ &\leq& \mathbb{E}\sup_{r\in[0,
s]}|u(r)-u^n(r)|^p_{H}\Big(\sum_{j=1}^{n-1}\int_{0}^{s}
r^{-2\nu}e^{-2r\lambda_{jn}^\frac\alpha2}dr\Big)^\frac p2. \nonumber \\
\end{eqnarray}
Arguing as in the proof of Lemma \ref{lem-HS-AN}, we get a constant
$K$ which depends only on  $ \nu $, such that
\begin{eqnarray}\label{est-final-second-Y-Yn}
\mathbb{E}\Big(\int_{0}^{s}(s&-&r)^{-2\nu}
\|E_ne^{-(s-r)\mathbb{A}^{\frac\alpha2}_{\mathcal{V}_{n-1}}}P_n\big(g(u(r))-\overbrace{g(u^n(r))}^n\big)\|^2_{HS}dr
\Big)^\frac p2 \Big) \nonumber \\
&\leq& K\mathbb{E}\sup_{r\in[0,
s]}|u(r)-u^n(r)|^p_{H}\Big(\sum_{j=1}^{\infty}j^{-\alpha(1-2\nu)}\Big)^\frac p2 \nonumber \\
&\leq& K\mathbb{E}\sup_{r\in[0,
s]}|u(r)-u^n(r)|^p_{H} \nonumber \\
\end{eqnarray}
since $ 0<\nu<\frac12-\frac1{2\alpha}$\del{series
$\sum_{j=1}^{\infty}j^{-\alpha(1-2\nu)}$ converges provided that $
$}. Now replacing \eqref{est-first-term-in-Y-Yn-final-corr} and
\eqref{est-final-second-Y-Yn} in \eqref{est-first-Y-Yn}, we obtain

\begin{eqnarray}\label{est-Y-Yn}
\mathbb{E}|Y(s)-Y_n(s)|^p_{H}&\leq& K_{p, T, b_\delta}\Big(
\Big(n^{-4\delta}+n^{-\alpha}\Big)^\frac p2+ \mathbb{E}\sup_{r\in[0,
s]}|u(r)-u^n(r)|^p_{H}\Big)\nonumber \\
\end{eqnarray}
consequently,
\begin{eqnarray}\label{est-final-Y-Yn}
\mathbb{E}|Y-Y_n|^p_{L^p(0,T; L^2(0,
1))}&=&\int_0^T\mathbb{E}|Y(s)-Y_n(s)|^p_{H}ds \nonumber \\ &\leq&
K_{p, T, b_\delta}\Big( \Big(n^{-4\delta}+n^{-\alpha}\Big)^\frac p2+
\int_0^T\mathbb{E}\sup_{r\in[0,
s]}|u(r)-u^n(r)|^p_{H}ds\Big)\nonumber \\
&\leq& K_{p, T, b_\delta}\Big( n^{-2\delta
p}+n^{-\frac{p\alpha}{2}}+ \int_0^T\mathbb{E}\sup_{r\in[0,
s]}|u(r)-u^n(r)|^p_{H}ds\Big)\nonumber \\
\end{eqnarray}
and by replacing \eqref{Inq. Y(s) in Lp(0,T; L2) Ex. case} and
\eqref{est-final-Y-Yn} in \eqref{Eq-ReRew-Est-C(B)-2}, we obtain

\begin{eqnarray}\label{Eq-final-ESupBt}
\mathbb{E}\Big[\sup_{[0,T]}B(t)\Big] &\leq& K_{p, T,
b_\delta}\Big[n^{-2\delta p}+n^{-\frac{p\alpha}{2}}+
\int_0^T\mathbb{E}\sup_{r\in[0, s]}|u(r)-u^n(r)|^p_{H}ds\Big].\nonumber \\
\end{eqnarray}

Finally, we join the estimations in \eqref{Eq-final-ESupBt} and
\eqref{Eq-EsupA}, we obtain,

\begin{eqnarray}\label{Eq-est-Mn-step2}
\mathbb{E}\big[\sup_{[0, T]}|u(t)-u^n(t)|_H^p\big] &\leq&
K_{p,T,|u_0|_{D(A^{\eta})}, b_\delta}^{p}\Big(
\int_0^T\mathbb{E}\sup_{r\in[0,
s]}|u(r)-u^n(r)|^p_{H}ds \nonumber\\
&+&n^{-2\eta p}+ n^{-2\delta p}+n^{-\frac{p\alpha}{2}}\Big).\nonumber \\
\end{eqnarray}
Thanks to Gronwall Lemma,

\begin{eqnarray*}
\mathbb{E}\big[\sup_{[0, T]}|u(t)-u^n(t)|_H^p\big] &\leq& K_{p,
T,|u_0|_{D(A^{\eta})},
b_\delta}^{p}\Big(n^{-2\eta p}+ n^{-2\delta p}+n^{-\frac{p\alpha}{2}}\Big).\nonumber \\
\end{eqnarray*}
Now it is easy to see that from the conditions on $ \eta $, that $
\frac \alpha2 \leq 2\eta$, hence

\begin{eqnarray}\label{Eq-est-Mn-step1}
\mathbb{E}\big[\sup_{[0, T]}|u(t)-u^n(t)|_H^p\big] &\leq&
K_{T,|u_0|_{D(A^{\eta})},
b_\delta}^{p}n^{-\xi p},\nonumber \\
\end{eqnarray}
where $ \xi$ is given by \eqref{Eq-xi-1}. Furthermore, if $ 1<\alpha
\leq 2$, then $ \frac\alpha2\leq \frac12+\frac\alpha4\leq
1<2\delta$, which implies that $ \xi = \frac\alpha2$.

\end{proof}


To prove Theorem \ref{Theom.Approx-2}, we will use the same scheme
as before, but with different estimations:

\begin{lem}\label{lem-estim-Salpha-discr-sgp-second}
For $ \gamma
>\frac1\alpha$, $ 0<\beta < 1$, there exists $ K>0$, such that $ \forall t \in (0, T]$,
\begin{eqnarray}\label{Eq-est-semiGroup-Salpha-Salphan-gamma-beta}
|S_\alpha(t)
-E_ne^{-t\mathbb{A}_{\CV_{n-1}}^{\frac\alpha2}}P_n|_{\mathcal{L}(H)}&+&\|S_\alpha(t)-
E_ne^{-t\mathbb{A}^{\frac\alpha2}_{\mathcal{V}_{n-1}}}P_n\|_{HS}\nonumber
\\ &\leq & K \big(n^{-\alpha\frac\gamma2}t^{-\frac\gamma2}+
n^{-\frac{\alpha-\alpha\beta}{2}}t^{-\frac{1+\alpha-\alpha\beta}{2\alpha}}\big),
\end{eqnarray}
\del{where
\begin{equation}\label{Eq-Def-Phi(t, n)-gamma-beta}
\phi_{\alpha, \delta, \gamma, \beta}(t, n) := \Big\{
\begin{array}{lr}
n^{-\alpha\frac\gamma2}t^{-\frac\gamma2}+
n^{-\frac{\alpha-\alpha\beta}{2}}t^{-\frac{1+\alpha-\alpha\beta-4\delta}{2\alpha}}, \, \, \, \, \, \,0\leq \delta\leq \frac14\\
n^{-2\delta}+
n^{-\frac{\alpha-\alpha\beta}{2}}t^{-\frac{1+\alpha-\alpha\beta-4\delta}{2\alpha}},
\, \, \, \, \, \, \, \, \, \, \frac14<\delta <
\frac14+\frac34\alpha.
\end{array}
\end{equation}}

\del{\begin{equation}\label{Eq-Def-Phi(t, n)} \phi_{\alpha, \eta}(t,
n) := \Big\{
\begin{array}{rl}
n^{-\alpha\frac\gamma2}t^{-\frac\gamma2}+
n^{-\frac{\alpha-\alpha\beta}{2}}t^{-\frac{1+\alpha-\beta\alpha}{2\alpha}}, \, \,\eta=0\\
n^{-2\eta}+
n^{-\frac{\alpha-\alpha\beta}{2}}t^{-\frac{1+\alpha-\beta\alpha-4\eta}{2\alpha}},
\, \, \eta\neq 0.
\end{array}
\end{equation}}
\del{\begin{equation}\label{Eq-Def-Phi(t, n)} \phi_{\alpha,
\delta}(t, n) := n^{-2\delta}+
n^{-\frac{\alpha}{2}}t^{-\frac{1+\alpha-4\delta}{2\alpha}}.
\end{equation}}
\end{lem}

\begin{proof}
The proof follows the same steps as in Lemmata
\ref{lem-estim-Salpha-discr-sgp} and \ref{lem-baic-est-G} and by
considering the following key estimation:
\begin{eqnarray}\label{Eq-Galpha-Galpha-n-partial2-second}
\sum_{j=1}^{n-1}|e^{-\lambda_j^\frac\alpha2t}-e^{-\lambda_{jn}^\frac\alpha2t}|^2&\leq&
K t^2n^{-4\alpha}\sum_{j=1}^{n-1}j^{6\alpha}
e^{-2^{\alpha+1}j^\alpha t}\nonumber\\
&\leq& K t^2n^{-\alpha+
\alpha\beta}\sum_{j=1}^{n-1}j^{3\alpha-\alpha\beta}
e^{-2^{\alpha+1}j^\alpha t}\nonumber\\
&\leq& K t^2n^{-\alpha+  \alpha\beta}\int_{0}^\infty
x^{3\alpha-\alpha\beta}
e^{-2^{\alpha+1}x^\alpha t}dx\nonumber\\
&\leq& K t^{-1-\frac1\alpha+\beta}n^{-\alpha+
\alpha\beta}\int_{0}^\infty y^{3\alpha-\alpha\beta}
e^{-2^{\alpha+1}y^\alpha }dy\nonumber\\
&\leq& K t^{-1-\frac1\alpha+\beta}n^{-\alpha+  \alpha\beta}.\nonumber\\
\end{eqnarray}
\end{proof}

\begin{lem}\label{lem G gamma-second}Provided $ p>\frac{2\alpha}{2\alpha-1}$, $ \frac1\alpha<\gamma
< 2\frac{p-1}{p}$ and $ \frac2p+\frac1\alpha-1<\beta $, there exists
$\nu$ satisfying:
$\max\{\frac{1+\alpha-\alpha\beta}{2\alpha}+\frac1p,
\frac\gamma2+\frac1p\}<\nu<1$, such that the operator $ G_{n, \nu}:
L^{p}(0, T; L^{2}(0, 1)) \rightarrow C([0, T]; L^{2}))$ given by
$$
G_{n, \nu}h(t)= \int_{0}^{t}(t-s)^{\nu-1}\big[S_{\alpha}(t-s)-
E_ne^{-(t-s)\mathbb{A}^{\frac\alpha2}_{\mathcal{V}_{n-1}}}P_n\big]h(s)\,ds,
\; h \in L^{p}(0, T; L^{2}(0, 1))
$$
is well defined, linear and bounded. Moreover, there exists a
constant $K_{p,\nu}>0$ such that for all $h \in L^{p}(0, T; L^{2}(0,
1))$
\begin{eqnarray}\label{ineq:Rg}
|G_{n, \nu}h|_{C([0, T]; L^2)}&\leq& K_{p,\nu}
(n^{-\alpha\frac\gamma2}+
n^{-\frac{\alpha-\alpha\beta}{2}})|h|_{L^{p}(0, T; L^{2}(0, 1))}.
\end{eqnarray}
\end{lem}

\begin{proof} Let $h \in L^{p}(0, T; L^{2}(0, 1))$. Then for
$t \in (0,T)$ we have
\begin{eqnarray*}
|G_{n, \nu}h(t)|_{L^2} &\leq& \int_{0}^{t}(t-s)^{\nu-1}\Vert
\big[S_{\alpha}(t-s)-
E_ne^{-(t-s)\mathbb{A}^{\frac\alpha2}_{\mathcal{V}_{n-1}}}P_n\big]\Vert_{\mathcal{L}(H)} |h(s)|_{L^2}\,ds.\nonumber\\
\end{eqnarray*}
From Lemma \ref{lem-estim-Salpha-discr-sgp-second} and by applying
the H\"older inequality, we get

\begin{eqnarray*} |G_{n, \nu}h(t)|_{L^2}&\leq&
K\int_{0}^{t}(t-s)^{\nu-1}
\big(n^{-\alpha\frac\gamma2}(t-s)^{-\frac\gamma2}+
n^{-\frac{\alpha-\alpha\beta}{2}}(t-s)^{-\frac{1+\alpha-\alpha\beta}{2\alpha}}\big)|h(s)|_{L^2}\,ds\nonumber\\
&\leq&
K\big(\int_{0}^{T}\big(s^{\nu-1}(n^{-\alpha\frac\gamma2}s^{-\frac\gamma2}+
n^{-\frac{\alpha-\alpha\beta}{2}}s^{-\frac{1+\alpha-\alpha\beta}{2\alpha}})\big)^{\frac{p}{p-1}}\,ds\big)^{\frac{p-1}{p}}
\big(\int_{0}^{T}|h(s)|_{L^2}^{p}\,ds\big)^{\frac{1}{p}}.\nonumber\\
\end{eqnarray*}
Thanks to the basic inequality: $ (x+y)^\theta\leq c_\theta
(x^\theta+ y^\theta)$, for $ \theta >1$ and $ x, y\geq 0,$ we have

\begin{eqnarray}\label{Conv-Int-sPower-second}
\int_{0}^{T}(s^{\nu-1}(n^{-\alpha\frac\gamma2}s^{-\frac\gamma2}+
n^{-\frac{\alpha-\alpha\beta}{2}}s^{-\frac{1+\alpha-\alpha\beta}{2\alpha}}))^{\frac{p}{p-1}}\,ds
&\leq &c_p(n^{-\alpha\frac\gamma2\frac{p}{p-1}}\int_{0}^{T}(s^{(\nu-1-\frac\gamma2)\frac{p}{p-1}}ds \nonumber\\
&+& n^{-\frac{\alpha-\alpha\beta}{2}\frac{p}{p-1}}\int_0^T
s^{(\nu-1-\frac{1+\alpha-\alpha\beta}{2\alpha})\frac{p}{p-1}}\,ds).
\nonumber\\
\end{eqnarray}
For $ p>\frac{2\alpha}{2\alpha-1}$, $ \frac1\alpha<\gamma <
2\frac{p-1}{p}$ there exists $\nu> \frac\gamma2+\frac1p$, such that
the first integral in the RHS of (\ref{Conv-Int-sPower-second})
converges. Furthermore, if  $ \beta > \frac2p+\frac1\alpha-1$ then
we can chose $\nu$ satisfying also the inequality  $\nu>
\frac{1+\alpha-\alpha\beta}{2\alpha}+\frac1p$ and hence the second
integral in the RHS of (\ref{Conv-Int-sPower-second}) converges
also. Hence for $\nu>
\max\{\frac{1+\alpha-\alpha\beta}{2\alpha}+\frac1p,
\frac\gamma2+\frac1p\}$, we have:
\begin{eqnarray*}
\big(\int_{0}^{T}(s^{\nu-1}(n^{-\alpha\frac\gamma2}s^{-\frac\gamma2}+
n^{-\frac{\alpha-\alpha\beta}{2}}s^{-\frac{1+\alpha-\alpha\beta}{2\alpha}}))^{\frac{p}{p-1}}\,ds\big)^{\frac{p-1}{p}}\leq
c_p(n^{-\alpha\frac\gamma2}+ n^{-\frac{\alpha-\alpha\beta}{2}}).
\nonumber\\
\end{eqnarray*}
Hence for all $ t \geq 0 $
\begin{eqnarray*}
|G_{\nu}h(t)|_{L^2}&\leq& K_{T, p}(n^{-\alpha\frac\gamma2}+
n^{-\frac{\alpha-\alpha\beta}{2}})|h|_{L^{p}(0, T; L^{2}(0, 1))}.
\nonumber
\end{eqnarray*}
\end{proof}

\subsection*{Proof of Theorem \ref{Theom.Approx-2}}
We arguing as in the proof of Theorem \ref{Theom.Approx}. We define

\begin{equation}\label{eq-Mn-2}
M_n:= ||u - u^n||^p_{T, p} = \mathbb{E}\sup_{t \in [0,
T]}|u(t)-u^n(t)|_{H}^{p}\leq c_p\Big(\mathbb{E}\sup_{t \in [0,
T]}A(t)+ \mathbb{E}\sup_{t \in [0, T]}B(t)\Big),
\end{equation}
where $ A(t)$ and $ B(t)$ are given respectively by
\eqref{eq-A(t)-1} and \eqref{eq-B(t)-1}. The term $\mathbb{E}\sup_{t
\in [0, T]}A(t)$ is estimated thanks to the inequality
\eqref{Eq-EsupA}. To estimate the term $B(t)$, we use Lemmata
\ref{lem R gamma} and \ref{lem G gamma-second}
\begin{eqnarray}\label{Eq-ReRew-Est-C(B)-second}
\mathbb{E}\Big[\sup_{[0,T]}B(t)\Big]&\leq & c_p\mathbb{E}\sup_{[0,T]}
\Big(|\int_0^t(t-s)^{\nu-1}S_\alpha(t-s)\big(Y(s)-Y_n(s)\big)ds|_{H}^p\nonumber \\
&+&|\int_0^t(t-s)^{\nu-1}\Big[S_\alpha(t-s)-
E_ne^{-(t-s)\mathbb{A}^{\frac\alpha2}_{\mathcal{V}_{n-1}}}P_n\Big]
Y_n(s)ds|_{H}^p\Big)\nonumber \\
&\leq&    K_{T, p}\Big[\mathbb{E}|Y-Y_n|^p_{L^{p}(0, T; L^{2}(0,
1))}+ (n^{-\frac{\alpha\gamma}2}+
n^{-\frac{\alpha-\alpha\beta}{2}})\mathbb{E}|Y_n|^p_{L^{p}(0, T;
L^{2}(0,
1))}\Big]\nonumber \\
\end{eqnarray}
Now we calculate $ \mathbb{E}|Y-Y_n|^p_{L^{p}(0, T; L^{2}(0, 1))}$
and $\mathbb{E}|Y_n|^p_{L^{p}(0, T; L^{2}(0, 1))}$. By the
Burkholder's inequality and  Lemma \ref{lem-HS-AN}, there exists a
constant $C_p>0$ such that
\begin{eqnarray}\label{Inq. Y(s) in Lp(0,T; L2) Ex. case-second}
\int_{0}^{T}\!\!\!\mathbb{E}\big|Y_n(s)
\big|_{L^2}^{p}\,ds \!\!\!&\leq &\!\!\!
K_p\int_{0}^{T}\!\!\!\mathbb{E}\Big(\int_{0}^{s}(s-r)^{-2\nu}\Vert
E_ne^{-(s-r)\mathbb{A}^{\frac\alpha2}_{\mathcal{V}_{n-1}}}P_ng(u(r))\Vert
_{HS}^{2}dr
\Big)^{\frac{p}{2}}\,ds \nonumber\\
\!\!\!&\leq &\!\!\!
K_p\int_{0}^{T}\!\!\!\mathbb{E}\Big(\int_{0}^{s}(s-r)^{-2\nu}\Vert
E_ne^{-(s-r)\mathbb{A}^{\frac\alpha2}_{\mathcal{V}_{n-1}}}P_n\Vert_{HS}^{2}\vert
g\vert_{L^\infty}^2dr
\Big)^{\frac{p}{2}}\,ds \nonumber\\
\!\!\! &\leq &\!\!\!
K_pb_0^{p}\int_{0}^{T}\!\!\!\Big(\int_{0}^{s}r^{-2\nu}\Vert
E_ne^{-r\mathbb{A}^{\frac\alpha2}_{\mathcal{V}_{n-1}}}P_n\Vert
_{HS}^{2}dr
\Big)^{\frac{p}{2}}\,ds \\
\!\!\!&\leq&\!\!\!K_pb_0^{p}T\!\Big(\sum_{k=1}^{+\infty} k^{-\alpha
(1-2\nu)} \Big)^{\frac{p}{2}}.
\nonumber
\end{eqnarray}
Since  $ 0 < \nu < \frac{1}{2}- \frac{1}{2 \alpha} $ what implies
that $\alpha(1-2\nu)>1$ we infer that the last term is finite. In
the aim to get an estimation to  $ \mathbb{E}|Y-Y_n|^p_{L^{p}(0, T;
L^{2}(0, 1))}$, we use the inequality \eqref{est-first-Y-Yn}. Let us
remark that the estimation \eqref{est-final-second-Y-Yn} remains
true for the second term in the RHS of \eqref{est-first-Y-Yn}. Let
us now estimate the first term in RHS of this inequality. We have:

\begin{eqnarray}\label{est-first-term-in-Y-Yn-second}
\int_{0}^{s}(s&-&r)^{-2\nu}\|S_\alpha(s-r)- E_ne^{-(s-r)\mathbb{A}^{\frac\alpha2}_{\mathcal{V}_{n-1}}}P_ng(u(r))\|_{HS}^2dr \nonumber \\
&\leq& \int_{0}^{s}(s-r)^{-2\nu}\|S_\alpha(s-r)-
E_ne^{-(s-r)\mathbb{A}^{\frac\alpha2}_{\mathcal{V}_{n-1}}}P_n\|_{HS}^2
\|g(u(r))\|_{\mathcal{L}(H)}^2dr\nonumber \\
&\leq& b_0^2\int_{0}^{s}(s-r)^{-2\nu}\| S_\alpha(s-r)-
E_ne^{-(s-r)\mathbb{A}^{\frac\alpha2}_{\mathcal{V}_{n-1}}}P_n\|_{HS}^2
dr\nonumber \\
\end{eqnarray}

Thanks to Lemma \ref{lem-estim-Salpha-discr-sgp-second}, we have

\begin{eqnarray*}
\int_{0}^{s}(s&-&r)^{-2\nu}\| S_\alpha(s-r)-
E_ne^{-(s-r)\mathbb{A}^{\frac\alpha2}_{\mathcal{V}_{n-1}}}P_n\|_{HS}^2
dr\nonumber \\
&\leq&K\Big(n^{-\alpha+\alpha\beta}\int_{0}^{s}r^{-2\nu-1-\frac1\alpha+\beta}dr+n^{-\alpha\gamma}\int_{0}^{s}r^{-2\nu-\gamma}dr\Big).\nonumber \\
\end{eqnarray*}
The last two integrals converge provided $ \nu
<\min\{\frac\beta2-\frac1{2\alpha}, \frac{1-\gamma}{2}\}$. Hence,
\begin{eqnarray}\label{est-first-term-in-Y-Yn-final-second}
\mathbb{E}\Big(\int_{0}^{s}(s&-&r)^{-2\nu}\|\big( S_\alpha(s-r)-
E_ne^{-(s-r)\mathbb{A}^{\frac\alpha2}_{\mathcal{V}_{n-1}}}
P_n\big)g(u(r))\|_{HS}^2dr \Big)^\frac p2\nonumber \\
&\leq&
b_0^pK_T\Big(n^{-\alpha+\alpha\beta}+n^{-\alpha\gamma}\Big)^\frac p2.\nonumber\\
\end{eqnarray}

By accumulating the conditions, the parameter $\nu$ should satisfy:
\begin{eqnarray}\label{eq-cond-on-nu-second}
\nonumber\max\{\frac{1+\alpha-\alpha\beta}{2\alpha}+\frac1p,
\frac\gamma2+\frac1p\}<\nu <\min\{\frac\beta2-\frac1{2\alpha},
\frac{1-\gamma}{2}, \frac12-\frac1{2\alpha}\}.\\
\end{eqnarray}

\del{Let us recall that the existence of $\nu$ satisfying the
inequalities \eqref{eq-cond-on-nu-second}, is thanks to the choice
of $ \alpha>2$, $ p>\frac{2\alpha}{\alpha-2}$, $
\frac1\alpha<\gamma<\frac12-\frac1p$ and $
\frac1\alpha+\frac1p+\frac12<\beta <1$. Let us remark that $
\frac2\alpha+\frac1p-1<\frac2\alpha+\frac2p<\frac1\alpha+\frac1p+\frac12<\beta
<1$ and $ \gamma
+\frac1\alpha+\frac2p<\frac1\alpha+\frac1p+\frac12<\beta.$}

Since $ \alpha>2$ and $ p>\frac{2\alpha}{\alpha-2}$, it is possible
to choose   $ \gamma \in (\frac1\alpha, \frac12-\frac1p)$ and $
\beta\in (\frac1\alpha+\frac1p+\frac12, 1)$ (notice that
$p>\frac{2\alpha}{\alpha-2}\Leftrightarrow
\frac1\alpha<\frac12-\frac1p\Leftrightarrow\frac1\alpha+\frac12+\frac1p<1$).
With such choice of $ \gamma $ and $ \beta$, we can find $\nu$
satisfying \eqref{eq-cond-on-nu-second}. In fact $
\frac2\alpha+\frac1p-1<\frac2\alpha+\frac2p<\frac1\alpha+\frac1p+\frac12<\beta
<1$ and $ \gamma
+\frac1\alpha+\frac2p<\frac1\alpha+\frac1p+\frac12<\beta.$  From
\eqref{est-final-second-Y-Yn} and
\eqref{est-first-term-in-Y-Yn-final-second}, we have
\del{\begin{eqnarray}\label{est-final-Y-Yn-second}
\mathbb{E}|Y-Y_n|^p_{L^p(0,T; L^2(0, 1))} &\leq& C_{p, T, b_0}\Big(
n^{-\alpha\frac p2+\alpha\beta\frac p2}+n^{-\alpha\gamma\frac p2}+
\int_0^T\mathbb{E}\sup_{r\in[0,
s]}|u(r)-u^n(r)|^p_{H}ds\Big)\nonumber \\
\end{eqnarray}
and}

\begin{eqnarray}\label{Eq-final-ESupBt-second}
\mathbb{E}\Big[\sup_{[0,T]}B(t)\Big] &\leq& K_{T, p,
b_0}\Big[n^{-\alpha\frac p2+\alpha\beta\frac
p2}+n^{-\alpha\gamma\frac p2}+
\int_0^T\mathbb{E}\sup_{r\in[0, s]}|u(r)-u^n(r)|^p_{H}ds\Big].\nonumber \\
\end{eqnarray}
Finally, we join the estimations in \eqref{Eq-final-ESupBt-second}
and \eqref{Eq-EsupA} and applying the Gronwall Lemma, we get
\begin{eqnarray}\label{Eq-est-Mn-step1-second}
||u - u^n||_{T, p}&:=&\Big(\mathbb{E}\big[\sup_{[0,
T]}|u(t)-u^n(t)|_H^p\big]\Big)^\frac1p \nonumber \\
&\leq&
K_{T,|u_0|_{D(A^{\eta})}, b_0}\Big(n^{-2p\eta}+n^{-\alpha\frac
p2}+n^{-\alpha\frac p2+\alpha\beta\frac
p2}+n^{-\alpha\gamma\frac p2}\Big)^\frac1p.\nonumber \\
\end{eqnarray}
To get a good precision of estimation, we push $ \beta $ to its
lower bound and $\gamma$ to its upper bound, i.e. we take $ \beta:=
\frac12+\frac1\alpha+\frac1p$ and $ \gamma =\frac12-\frac1p$, then
we get:
\begin{eqnarray}\label{Eq-est-Mn-step1-second-pure}
||u - u^n||_{T, p} &\leq& K_{T,|u_0|_{D(A^{\eta})},
b_0}\Big(n^{-2p\eta}+n^{-\alpha\frac p2}+n^{-\alpha\frac
p2+\alpha\frac p2(\frac12+\frac1p+\frac1\alpha)}+n^{-\alpha\frac
p2(\frac12-\frac1p)}\Big)^\frac1p\nonumber \\
&\leq& K_{T,|u_0|_{D(A^{\eta})},
b_0}\Big(n^{-2p\eta}+n^{-(\frac\alpha p4-\frac\alpha2-\frac
p2)}\Big)^\frac1p\nonumber \\
\end{eqnarray}
Thanks to the condition \eqref{eq-cond-eta}, we have $ \frac{\alpha
p}4-\frac\alpha2-\frac p2<\frac p2+\frac{\alpha p}2<2\eta p$. Now,
it is easy to get the inequality
\eqref{Eq-mean-rate-estimation-second}.

\appendix

\del{\section*{Definitions}
\begin{defn} Let us  recall that the Hilbert-Schmidt norm of a
linear operator $\Lambda$ in a Hilbert space $H$ is by definition
equal to $$\Vert\Lambda \Vert _{HS}=
\Big(\sum_{k=1}^{\infty}\Vert\Lambda e_k\Vert
_{L^2}^2\Big)^\frac12,$$ where $\{e_k\}_{k}$ is an ONB of $H$. The
RHS of the equality above is in fact independent of the choice of
the ONB. A linear operator $\Lambda$ in a Hilbert space $H$ is
called to be a Hilbert-Schmidt operator if and only if $\Vert\Lambda
\Vert _{HS}$ is finite.
\end{defn}}

\section{Proof of Lemma \ref{appendix-E-n-Stochas}}\label{Appendix-proof-Lemma}
We apply on the both sides of equation \eqref{Eq-Solu-integ} the
operator $E_n$, we get
\begin{equation}\label{Eq-Solu-integ-App}
E_nu_n(t)= E_n
\big(e^{-t\mathbb{A}_{\CV_{n-1}}^{\frac\alpha2}}u_n(0)\big)+
E_n\Big(\int_{0}^{t}e^{-(t-s)\mathbb{A}_{\CV_{n-1}}^{\frac\alpha2}}g_n(u_n(s))dW_n(s)\Big).
\end{equation}
Using the definitions of the two operators $ P_n$ and $ E_n$ and by
the fact that $ P_n E_n= I$, the RHS in \eqref{Eq-Solu-integ-App} is
equal to $ u^n(t)$. Moreover
$$E_n e^{-t\mathbb{A}_{\CV_{n-1}}^{\frac\alpha2}}u_n(0)= E_ne^{-t\mathbb{A}_{\CV_{n-1}}^{\frac\alpha2}
}P_nu_0.$$ Let us now explain here how to get the stochastic term.
We denote by $(\cdot)_k$  the component of  a vector and by
$(\cdot)_{kj}$ the component of  a matrix. We have
\begin{eqnarray*}
E_n\Big(\int_{0}^{t}e^{-(t-s)\mathbb{A}_{\CV_{n-1}}^{\frac\alpha2}}&g_n(u_n(s))&dW_n(s)\Big):=
\sum_{k=1}^{n-1}\Big(\int_{0}^{t}e^{-(t-s)\mathbb{A}_{\CV_{n-1}}^{\frac\alpha2}}g_n(u_n(s))dW_n(s)\Big)_ke_k \nonumber\\
&=&
\sum_{k=1}^{n-1}\int_{0}^{t}\sum_{j=1}^{n-1}\Big(e^{-(t-s)\mathbb{A}_{\CV_{n-1}}^{\frac\alpha2}}g_n(u_n(s))\Big)_{kj}dB_j(s)e_k\nonumber\\
&=&
\sum_{j=1}^{n-1}\int_{0}^{t}\sum_{k=1}^{n-1}\big(e^{-(t-s)\mathbb{A}_{\CV_{n-1}}^{\frac\alpha2}}g_n(u_n(s))\big)_{kj}e_kdB_j(s).\nonumber\\
\end{eqnarray*}
But $$
\sum_{k=1}^{n-1}\big(e^{-(t-s)\mathbb{A}_{\CV_{n-1}}^{\frac\alpha2}}g_n(u_n(s))\big)_{kj}e_k
= E_n
\big(e^{-(t-s)\mathbb{A}_{\CV_{n-1}}^{\frac\alpha2}}g_n(u_n(s))\big)_{j},$$
where
$\big(e^{-(t-s)\mathbb{A}_{\CV_{n-1}}^{\frac\alpha2}}g_n(u_n(s))\big)_{j}$
is the column $"j"$ of the matrix
$e^{-(t-s)\mathbb{A}_{\CV_{n-1}}^{\frac\alpha2}}g_n(u_n(s))$.

Hence we have the first result:
\begin{eqnarray*}
E_n\Big(\int_{0}^{t}e^{-(t-s)\mathbb{A}_{\CV_{n-1}}^{\frac\alpha2}}g_n(u_n(s))dW_n(s)\Big)=
\sum_{j=1}^{n-1}\int_0^t
E_n(e^{-(t-s)\mathbb{A}_{\CV_{n-1}}^{\frac\alpha2}}g_n(u_n(s)))_jdB_j(s)\\.
\end{eqnarray*}

We know from the basic calculus on matrices that:
$$\big(e^{-(t-s)\mathbb{A}_{\CV_{n-1}}^{\frac\alpha2}}g_n(u_n(s))\big)_{j}= e^{-(t-s)\mathbb{A}_{\CV_{n-1}}^{\frac\alpha2}}
\big(g_n(u_n(s))\big)_{j}. $$ By the definition of the matrix $g_n$,
the  column $"j"$ of the matrix $g_n(u_n(s))$: ($ (g_n(u_n(s)))_j$)
is equal to $P_n((g\circ E_n)e_j)$, hence we have

\begin{eqnarray*}
E_n\Big(\int_{0}^{t}e^{-(t-s)\mathbb{A}_{\CV_{n-1}}^{\frac\alpha2}}g_n(u_n(s))dW_n(s)\Big)&=&
\sum_{j=1}^{n-1}\int_0^t
E_n(e^{-\mathbb{A}_{\CV_{n-1}}^{\frac\alpha2}(t-s)}g_n(u_n(s)))_jdB_j(s)\\
&=&\sum_{j=1}^{n-1}\int_{0}^{t}E_ne^{-\mathbb{A}_{\CV_{n-1}}^{\frac\alpha2}(t-s)}P_ng(u^n(s))e_jdB_j(s)
\end{eqnarray*}
and we denote these integrals by

\begin{eqnarray*}
\int_{0}^{t}E_ne^{-\mathbb{A}_{\CV_{n-1}}^{\frac\alpha2}(t-s)}P_ng(u^n(s))dW^n(s).
\end{eqnarray*}
To get the estimation of $u^n(t)$, we  use Lemma
\ref{lem-properties-Pn-En} and Theorem \ref{Theom.Exist.approx
Solu.}, then we have
$$\mathbb{E} \sup_{t\in [0,
T]}|u^{n}(t)|_{L^2}^{p} \leq \mathbb{E} \sup_{t\in [0,
T]}||E_n|||u_{n}(t)|_{\mathbb{R}^{n-1}}^{p} \leq C_{T, n,
||g||}(1+\mathbb{E}|u_0|_{L^{2}}^{p}), \text{ for each }T>0.$$

\section{Proof of Lemma \ref{lem-Cumm}}
In fact, the operators $S_\alpha(t)$,
$E_ne^{-\mathbb{A}_{\CV_{n-1}}^{\frac\alpha2}t}P_n, $ and
$A^{-\delta} $ are bounded and we have $A^{-\delta} S_\alpha(t) e_j=
S_\alpha(t)A^{-\delta}e_j$ and $
E_ne^{-\mathbb{A}_{\CV_{n-1}}^{\frac\alpha2}t}P_n A^{-\delta}e_j =
A^{-\delta}E_ne^{-\mathbb{A}_{\CV_{n-1}}^{\frac\alpha2}t}P_ne_j$,
for all $ e_j, j\in \mathbb{N} $.

\section{Lemma \ref{lem-eigenvalues-fract-and-initial-ops}}
\begin{lem}\label{lem-eigenvalues-fract-and-initial-ops}
Let $(\lambda_k)_{k\geq0}$ be the sequence of  the eigenvalues
corresponding to the eigenfunctions $(e_k)_{k\geq0}$ of the positive
operator $A$. Then $(e_k)_{k\geq0}$ are also eigenfunctions of $
A^{\frac\alpha2}$ corresponding to the eigenvalues $
(\lambda^{\frac\alpha2}_k)_{k\geq0}$.
\end{lem}

\begin{proof}
Using the definition of the fractional operator
\eqref{Eq-Fract-Pazy}:
\begin{eqnarray}
A^{\frac\alpha2}e_k&=&
\frac{\sin\frac{\alpha\pi}2}{\pi}\int_0^\infty
t^{\frac\alpha2-1}A(tI+A)^{-1}e_k dt\nonumber \\
&=& \frac{\sin\frac{\alpha\pi}2}{\pi}\int_0^\infty
t^{\frac\alpha2-1}\lambda_k(t+\lambda_k)^{-1}dt e_k \nonumber \\
&=&
\lambda_k^{\frac\alpha2}\big(\frac{\sin\frac{\alpha\pi}2}{\pi}\int_0^\infty
\xi^{\frac\alpha2-1}(t+\xi)^{-1}d\xi \big) e_k. \nonumber \\
\end{eqnarray}
By the residues theory, we have
$$(1-e^{i\pi\alpha})\int_0^\infty
\xi^{\frac\alpha2-1}(t+\xi)^{-1}d\xi = 2\pi iRes(-1,
\xi^{\frac\alpha2-1}(t+\xi)^{-1})=-2\pi ie^{i\pi\frac\alpha2}$$
so,
$$\int_0^\infty
\xi^{\frac\alpha2-1}(t+\xi)^{-1}d\xi =
\frac{\pi}{\sin\frac{\alpha\pi}2}.$$

\end{proof}


\begin{thebibliography}{30}

\bibitem{GongyAlabert-06}  Alabert A. and Gy\"{o}ngy I. \textit{On numerical approximation of stochastic {B}urgers'
equation}, From stochastic calculus to mathematical finance,
Springer Berlin 1--15 (2006).

\bibitem{Bilerandal} Biler P., Funaki, T. and  Woyczynski W. A.
\textit{Fractal {B}urgers' equations.} J. {D}ifferential {E}quations
148, 9--46 (1998).

\bibitem{Brenner-Scott} Brenner S. C. and Scott L. R. \textit{The mathematical theory of
finite element methods.} Texts in Applied Mathematics 15.
Springer-Verlag, New York (2002).

\bibitem{BrzezniakDebbi1} Brze{\'z}niak, Z. and Debbi L. \textit{On Stochastic {B}urgers
Equation Driven by a Fractional Power of the {L}aplacian and
space-time white noise. } Stochastic Differential Equation: Theory
and Applications, A volume in Honor of Professor Boris L. Rozovskii.
Edited by P. H. Baxendale and S. V. Lototsky, 135--167 (2007).

\bibitem{Caffarelli-2009} Caffarelli L.A. \textit{Some nonlinear problems involving
non-local diffusions.} ICIAM 07-6th Intern. Congress on Industrial
and Applied Math.,  Eur. Math. Soc. Zurich 43--56 (2009).


\bibitem{Caffarelli-Vasseur2010} Caffarelli L.A. and
Vasseur A. \textit {Drift diffusion equations with fractional
diffusion and the quasi-geostrophic equation.} Ann. of Math. 2, 171
no. 3, 1903-1930 (2010).

\bibitem{Anh-09} Chen S., Liu F., Zhuang P. and Anh V. \textit{Finite difference
approximations for the fractional Fokker-Planck equation.} Appl.
Math. Model. 33 no. 1, 256--273 (2009).

\bibitem{DaPZa-92} Da Prato G. and Zabczyk J. \textit{Stochastic Equations in Infinte
Dimensions.} Springer,  Combridge university press (1992).

\bibitem{DebbiDozzi1} Debbi L. and  Dozzi M. \textit{On The Solution of Non Linear
Stochastic Fractional Partial Differential Equations.} Stochastic
Process. Appl. 115 no 11, 1764--1781 (2005).

\bibitem{Droniou-09} Droniou J.\textit{A numerical method for fractal conservation laws.}, Math. Comp. 79  no.
269, 95-124 (2010).

\bibitem{Friedman-75} Friedman A. \textit{Stochastic differential equations and
applications.} Vol. 1. Probability and Mathematical Statistics Vol.
28. Academic Press, New York-London 1975.

\bibitem{Gyongy-Millet-09} Gy\"ongy I. and  Millet A. \textit{Rate of convergence of space
time approximations for stochastic evolution equations.}  Potential
Anal. 30  no. 1, 29-64 (2009).

\bibitem{Gyongy-Millet-05} Gy\"ongy I. and  Millet A. \textit{On discretization schemes for
stochastic evolution equations.}  Potential Anal.  23  no. 2,
99-134 (2005).

\bibitem{Hausenblas-03} Hausenblas E. \textit{Approximation for semilinear stochastic
evolution equations.} Potential Anal. 18 no. 2, 141--186 (2003).

\bibitem{LM-72-i} Lions J.-L. and Magenes E. \textit{Non-homogeneous boundary value problems
and applications. {V}ol. {I}} Springer-Verlag, New York 1972.


\bibitem{Meerschaert-06} Meerschaert M., Tadjeran C. and  Scheffler H. P. \textit{A
second-order accurate numerical approximation for the fractional
diffusion equation.} J. Comput. Phys. 213 no. 1, 205--213 (2006).

\bibitem{LeonZas-02} Leoncini X. and Zaslavsky G. M. \textit{Jets, stickiness and
anomalous transport.} Phys. Rev. E (3)  65 no. 4, 046216, 16 pp65,
(2002).

\bibitem{Mueller-98} Mueller C. \textit{The heat equation with L\'evy
noise.} Stoch. Proc. Appl. 74, 67--82 (1998).

\bibitem{Pablo-VazquezA-Arxiv} Pablo A. D., Quiros F., Rodriguez A. and VazquezA J. L. \textit{Fractional porous medium equation.}
arXiv:1001.2383v1 [math.AP] 14 Jan (2010).


\bibitem{Pazy} Pazy A.\textit{Semigroups of Linear Operators and Applications to Partial Differential Equations}, Springer-Verlag,
New York  1983.


\bibitem{Printems-01} Printems J.\textit{On the discretization in time of parabolic stochastic partial differential equations},
M2AN Math. Model. Numer. Anal. 35 no. 6, 1055--1078, (2001).

\bibitem{Schwab-10} Schneider R., Reichmann O. and Schwab C. \textit{Wavelet solution of variable order pseudodifferential equations.}
Calcolo 47 no 2, 65--101 (2010).

\bibitem{SW1} Schneider W. R. and Wyss W.\textit{Fractional Diffusion and Wave Equations.} J. Math. Phys, 30 no 1,
134--144 (1989).

\bibitem{Anh-Tuner-08} Shen S., Liu F., Anh V. and Turner I. \textit{The fundamental solution and numerical solution of the Riesz
fractional advection-dispersion equation.} IMA J. Appl. Math. 73 no
6, 850--872 (2008).

\bibitem{Woyczynski-Num-approx-fract-05} Stanescu D., Kim D. and  Woyczynski W. A. \textit{Numerical study of interacting particles
approximation for integro-differential equations.} J. Comput.
Phys.206 no. 2, 706-726 (2005).

\bibitem{Sug} Sugimoto N. \textit{Generalized {B}urgers Equations and fractional calculus.} Nonlinear {W}ave {M}otion (A. Jeffrey, {E}d)
162--179 (1989).

\bibitem{SugKak} Sugimoto N. and Kukatani T.\textit{Generalized {B}urgers Equations
for Nonlinear Viscoelastic Waves.} Wave Motion 7, 447--458 (1985).

\bibitem{Ta} Taylor M. E. \textit{Pseudodifferential Operators.} Princeton University Press
1981.

\bibitem{Triebel-95} Triebel H. \textit{Interpolation theory, function spaces, differential
operators.} Second Edition, Heidelberg 1995.

\bibitem{TrumanWu-06} Truman A. and Wu J.L.\textit{Fractal Burgers' equation driven by Lévy
noise.}  Stochastic partial differential equations and
applicationsVII,  295-310, Lect. Notes Pure Appl. Math. 245
Chapman \& Hall/CRC Boca Raton FL 2006.

\bibitem{Westphal-73} Westphal U. \textit{An aproach to fractional powers of operators
via fractional difference.} Proc. London Math. Soc. 29 no 3,
557--576 (1974).

\bibitem{ZasABD-95} Zaslavsky M. G. and Abdullaev S. S. \textit{Scaling Property and
Anomalous Transport of Particles Inside the Stochastic Layer.} Phys.
Rev. E 51 no. 5, 3901--3910 (1995).









\end{thebibliography}
\end{document}